\documentclass[11pt,letterpaper, ]{article}

\newcommand{\bea}{\begin{eqnarray*}}
\newcommand{\eea}{\end{eqnarray*}}
\newcommand{\be}{\begin{eqnarray}}
\newcommand{\ee}{\end{eqnarray}}
\newcommand{\ed}{\end{document}}
\newcommand{\no}{\noindent}

\newcommand{\btab}{\begin{tabular}}
\newcommand{\etab}{\end{tabular}}
\newcommand{\bc}{\begin{center}}
\newcommand{\ec}{\end{center}}

\newcommand{\bi}{\begin{itemize}}
\newcommand{\ei}{\end{itemize}}
\newcommand{\bfi}{\begin{figure}}
\newcommand{\efi}{\end{figure}}
\newcommand{\ben}{\begin{enumerate}}
\newcommand{\een}{\end{enumerate}}
\newcommand{\bdes}{\begin{description}}
\newcommand{\edes}{\end{description}}
\newcommand{\bay}{\begin{array}}
\newcommand{\eay}{\end{array}}

\def\bco{\iffalse}

\def\cp{\citep}

\def\i1{_{i1}}

\def\cB{{\mathcal B}}

\def\s2{\sum_{k=1}^\infty\lambda_k\phi_k(t) \phi_k(t)}
\def\s1{\sum_{k=1}^\infty\lambda_k\phi^{(1)}_k(s) \phi_k(t)}

\def\pt{\phi_k(t)}
\def\ps{\phi_k(s)}
\def\pt1{\phi^{(1)}_k(t)}
\def\ps1{\phi^{(1)}_k(s)}

\def\singlespace{\def\baselinestretch{1}\@normalsize}

\newcommand{\single}{\renewcommand{\baselinestretch}{1.2}\normalsize}
\newcommand{\double}{\renewcommand{\baselinestretch}{1.5}\normalsize}

\usepackage{geometry,float}
\usepackage{amssymb}
\usepackage{amstext}
\usepackage{graphics}
\usepackage{graphicx}
\usepackage{amsthm}
\usepackage{amsmath}
\usepackage{mathtools}
\usepackage{enumitem}
\usepackage{graphics}
\usepackage{rotating}
\usepackage[linesnumbered,lined,ruled]{algorithm2e}
\usepackage[colorlinks=true, allcolors=blue, draft=false]{hyperref}

\usepackage{natbib}

\AtBeginDocument{%
  \addtolength\abovedisplayskip{0\baselineskip}%
  \addtolength\belowdisplayskip{0\baselineskip}%
}

\bibpunct{(}{)}{;}{a}{}{,}
\graphicspath{{figures/}}

\newcommand*{\diff}[2]{\frac{\partial #1}{\partial #2}} 
\AtBeginDocument{}
\makeatletter
\newsavebox{\mybox}\newsavebox{\mysim}
\newcommand{\distas}[1]{%
  \savebox{\mybox}{\hbox{\kern3pt$\scriptstyle#1$\kern3pt}}%
  \savebox{\mysim}{\hbox{$\sim$}}%
  \mathbin{\overset{\text{#1}}{\kern\z@\resizebox{\wd\mybox}{\ht\mysim}{$\sim$}}}%
}
\makeatother

\newcommand*{\conv}[1]{\mathrel{\mathop{\longrightarrow}\limits^{#1}}}

\DeclarePairedDelimiterX{\inner}[2]{\langle}{\rangle}{#1, #2}
\DeclarePairedDelimiterX{\innerp}[2]{\langle}{\rangle_p}{#1, #2}
\newcommand{\innerR}[2]{{#1}^T{#2}}

\DeclareMathOperator*{\argmin}{arg\,min}

\newcommand{\ntoinf}{n \rightarrow \infty}
\newcommand{\deltatozero}{\delta \downarrow 0}
\newcommand{\toinf}{\rightarrow \infty}
\newcommand{\tozero}{\rightarrow 0}

\newcommand{\sumkinf}{\sum_{k=1}^\infty}
\newcommand{\sumjJ}{\sum_{j=1}^J}
\newcommand{\sumkK}{\sum_{k=1}^K}
\newcommand{\sumin}{\sum_{i=1}^n}

\newcommand{\bbR}{\mathbb{R}}

\newcommand{\bbN}{\mathbb{N}}

\newcommand{\bbH}{\mathbb{H}}

\newcommand{\tr}{\,\mathrm{tr}}

\newcommand{\FVE}{\,\mathrm{FVE}}

\newcommand{\cov}{\,\mathrm{cov}}

\DeclareMathOperator{\Exp}{Exp}
\DeclareMathOperator{\Log}{Log}

\newcommand{\inj}{\,\mathrm{inj}}

\renewcommand*{\vec}[1]{\mathbf{#1}}

\newcommand*{\bPsi}{\boldsymbol{\Psi}}
\newcommand*{\hbPhi}{\hat{\boldsymbol{\Phi}}}

\newcommand{\by}{\vec{y}}

\newcommand{\hmu}{\hat{\mu}}

\newcommand{\hV}{\hat{V}}
\newcommand{\hU}{\hat{U}}
\newcommand{\hK}{\hat{K}}

\newcommand{\hG}{\hat{G}}

\newcommand{\hxi}{\hat{\xi}}

\newcommand{\hlambda}{\hat{\lambda}}
\newcommand{\hphi}{\hat{\phi}}

\newcommand{\hX}{\hat{X}}
\newcommand{\hFVE}{\widehat{\FVE}}
\newcommand{\hnu}{\hat{\nu}}

\newcommand{\tG}{\tilde{G}}

\newcommand*{\cA}{\mathcal{A}}
\newcommand*{\cT}{\mathcal{T}}

\newcommand*{\cM}{\mathcal{M}}
\newcommand*{\cN}{\mathcal{N}}
\newcommand*{\cK}{\mathcal{K}}

\newcommand*{\cV}{\mathcal{V}}

\newcommand*{\cC}{\mathcal{C}}
\newcommand*{\cX}{\mathcal{X}}

\usepackage{amsthm}
\theoremstyle{compact.definition}
\newtheorem{Theorem}{Theorem}

\newtheorem{Corollary}{Corollary}
\newtheorem{Proposition}{Proposition}

\makeatletter
\def\th@newremark{\th@remark\thm@headfont{\bfseries}} 
\makeatletter
\theoremstyle{newremark}
\newtheorem{Remark}{Remark}

\theoremstyle{compact.remark}

\newcommand{\LFPCA}{$L^2$~FPCA{}}

\newcommand{\muM}{{\mu_\cM}}

\newcommand{\hmuM}{{\hmu_\cM}}

\newcommand{\muMt}{{\mu_\cM(t)}}

\newcommand{\dM}{{d_\cM}}

\newcommand{\dX}{{d_\cX}}

\newcommand{\muMs}{{\mu_\cM(s)}}
\newcommand{\Tmu}{T_{\mu_\cM(t)}}
\newcommand{\logmu}{\log_{\mu_\cM(t)}}
\newcommand{\loghmu}{\log_{\hmu_\cM(t)}}
\newcommand{\expmu}{\exp_{\mu_\cM(t)}}

\newcommand{\hmuMt}{{\hmu_\cM(t)}}

\newcommand{\hXiK}{\hX_{iK}}
\newcommand{\intcT}{\int_\cT}

\newcommand{\cMK}{\cM_K}
\newcommand{\cVK}{\cV_K}
\newcommand{\lambdamin}{\lambda_{\min}}
\newcommand{\inft}{\inf_{t \in \cT}}
\newcommand{\supt}{\sup_{t \in \cT}}

\newcommand{\suppK}{\sup_{p \in \cK}}
\newcommand{\supts}{\sup_{t,s \in \cT}}
\newcommand{\infdp}{\inf_{p:\, \dM(p, \muMt)>\epsilon}}
\newcommand{\infdpK}{\inf_{\substack{p\in\cK,\\\dM(p, \muMt)>\epsilon}}}

\newcommand{\supdpK}{\sup_{\substack{p\in\cK \\ \dM(p, \muMt)>\epsilon}}}

\newcommand{\suptsdelta}{\sup_{|t-s| < \delta}}
\newcommand{\suptspq}{\sup_{\substack{|t-s| < \delta \\ p,q\in\cK \\ \dM(p, q) < \delta}}}
\newcommand{\oneovern}{\frac{1}{n}}

\newcommand{\tinT}{t \in \cT}

\newcommand{\pinM}{p \in \cM}
\newcommand{\sqrtn}{\sqrt{n}}

\newcommand{\normF}[1]{\left\lVert#1\right\rVert_F}
\newcommand{\normR}[1]{\lVert#1\rVert_{E}}
\newcommand{\norm}[1]{\left\lVert#1\right\rVert}

\newcommand{\injmuMt}{\inj_{\muMt}}

\newcommand{\bbRdz}{\bbR^{d_0}}
\newcommand{\SON}{SO$(N)$}

\newcommand{\dtau}{d_{\tau}}

\newcommand{\tnu}{\tilde{\nu}}

\newcommand{\hbVi}{\hat{\mathbf{V}}_i}
\newcommand{\hbG}{\hat{\mathbf{G}}}

\double

\def\references{\bibliography{sphere}}

\begin{document}

\thispagestyle{empty}  \bc {\bf \sc \Large Principal Component Analysis for Functional Data on Riemannian Manifolds and Spheres} \vspace{0.25in}\ec

\vspace*{0.1in} \centerline{\today} \vspace*{0.1in}

\vspace*{0.1in} \centerline{Short title: Functional Data on Riemannian Manifolds}
\vspace*{0.1in}

\begin{center}

Xiongtao Dai\\

Department of Statistics \\ %
University of California, Davis \\ %
Davis, CA 95616 U.S.A. \\ %
Email: dai@ucdavis.edu

\vspace*{0.3in} 
Hans-Georg M\"uller

Department of Statistics \\ %
University of California, Davis \\ %
Davis, CA 95616 U.S.A. \\ %
Email: hgmueller@ucdavis.edu \\

\end{center}

\thispagestyle{empty}\vfill

 \vspace{-5.9cm}\noindent\rule{\textwidth}{0.5pt}\vspace{.01cm}
{\small Research supported by NSF grants 
DMS-1407852 and DMS-1712864.}

\vspace{-.275cm}

\noindent {\small We thank FlightAware for the permission to use the flight data in \autoref{s:data}.
}

\newpage \pagenumbering{arabic}

 \bc {\bf ABSTRACT} \ec

\no Functional data analysis on nonlinear manifolds has drawn recent interest. Sphere-valued functional data, which are encountered for example as movement trajectories on the surface of the earth, are an important special case. We consider an intrinsic principal component analysis for smooth Riemannian manifold-valued functional data and study its asymptotic properties. Riemannian functional principal component analysis (RFPCA) is carried out by first mapping the manifold-valued data through Riemannian logarithm maps to tangent spaces around the time-varying Fr\'echet mean function, and then performing a classical multivariate functional principal component analysis on the linear tangent spaces. Representations of the Riemannian manifold-valued functions and the eigenfunctions on the original manifold are then obtained with exponential maps. The tangent-space approximation through functional principal component analysis is shown to be well-behaved in terms of controlling the residual variation if the Riemannian manifold has nonnegative curvature. Specifically, we  derive a central limit theorem for the mean function, as well as root-$n$ uniform convergence rates for other model components, including the covariance function, eigenfunctions, and functional principal component scores. Our applications include  a novel framework for the analysis of longitudinal compositional data, achieved by mapping longitudinal compositional data to trajectories on the sphere, illustrated with longitudinal fruit fly behavior patterns. Riemannian functional principal component analysis is shown to be superior in terms of trajectory recovery in comparison to an unrestricted functional principal component analysis in applications and simulations and is also found to produce principal component scores that are better predictors for classification  compared to traditional functional functional principal component scores. 

\vspace{0.1in}

\no {\it Key words and phrases:} Compositional Data, Dimension Reduction, Functional Data Analysis, Functional Principal Component Analysis, Principal Geodesic Analysis, Riemannian Manifold, Trajectory, Central Limit Theorem, Uniform Convergence\\

\no {\bf MSC2010 Subject Classification:} Primary 62G05; secondary 62G20, 62G99.

\thispagestyle{empty} \vfill

\section{Introduction}
Methods for functional data analysis in a linear function space \citep{wang:16} or on a nonlinear submanifold \citep{lin:17} have been much studied in recent years. Growth curve data \citep{rams:05} are examples of functions in a linear space, while densities \citep{knei:01} and longitudinal shape profiles \citep{kent:01} lie on nonlinear manifolds. 
Since random functions usually lie in an intrinsically infinite dimensional linear or nonlinear space, dimension reduction techniques, in particular functional principal component analysis, play a central role in representing the random functions \citep{pete:16} and in other supervised/unsupervised learning tasks.  
Methods for analyzing non-functional data on manifolds have also been well developed over the years, such as data on spheres \citep{fish:87}, Kendall's shape spaces \citep{kend:09,huck:10}, and data on other classical Riemannian manifolds \citep{corn:17}; for a comprehensive overview of nonparametric methods for data on manifolds see \cite{patr:15}. Specifically,  versions of principal component analysis methods that adapt to the Riemannian or spherical geometry, such as principal geodesic analysis \citep{flet:04} or nested spheres \citep{huck:16},  have substantially advanced the study of data on manifolds. 

However, there is much less known about functional data, i.e., samples of random trajectories, that assume values  on manifolds, even though such data are quite common. An example is \cite{tels:16}, who considered the extrinsic mean function and warping for functional data lying on SO(3). Examples of data lying on a Euclidean sphere include geographical data \citep{zhen:15} on $S^2$, directional data on $S^1$ \citep{mard:09}, and square-root compositional data \citep{huck:16}, for which we will study longitudinal/functional versions in \autoref{s:comp}.
Sphere-valued functional data naturally arise when data on a sphere have  a time component, such as in recordings of airplane flight paths or animal migration trajectories. 
Our main goal is to extend and study the dimension reduction that is afforded by the popular functional principal component analysis (FPCA) in Euclidean spaces to the case of samples of  smooth curves that  lie on a smooth Riemannian manifold,  taking into account the underlying geometry. 
 
Specifically,  Riemannian Functional Principal Component Analysis (RFPCA) is shown to serve as an  intrinsic principal component analysis of Riemannian manifold-valued functional data. Our approach provides a theoretical framework and differs from existing methods for functional data analysis that involve manifolds, e.g., a proposed smooth principal component analysis for functions whose domain is on a two-dimensional manifold, motivated by signals on the cerebral cortex \cp{lila:16},  nonlinear manifold representation of $L^2$ random functions themselves lying on a low-dimensional but unknown manifold \cp{mull:12:1}, or functional predictors lying on a smooth low-dimensional manifold \cp{lin:17}. While there have been closely related computing and  application oriented proposals, including functional principal components on manifolds in discrete time, 
a systematic approach and theoretical analysis within a statistical modeling framework does not exist yet, to the knowledge of the authors.  
Specifically, in the engineering literature, dimension reduction for Riemannian manifold-valued motion data has been considered \citep{rahm:05,tour:09,anir:15}, where for example in the latter paper the time axis is  discretized, followed by multivariate dimension reduction techniques such as principal component analysis on the logarithm mapped data; these works emphasize specific applications and do not provide  theoretical justifications.   
The basic challenge  is to adapt inherently  linear methods such as functional principal component analysis (FPCA) to curved spaces.

RFPCA is an approach intrinsic to a given smooth Riemannian manifold and proceeds through time-varying geodesic submanifolds on the given manifold by minimizing total residual variation as measured by geodesic distance on the given manifold.   Since the mean of manifold-valued functions in the $L^2$ sense is usually extrinsic, i.e., does not lie itself  on the manifold in general, for an intrinsic analysis the mean function needs to be carefully defined, for which we adopt the intrinsic Fr\'echet mean, assuming that 
it is uniquely determined. RFPCA is implemented by first mapping the manifold valued trajectories that constitute the  functional data onto the linear tangent spaces using logarithm maps around the mean curve at a current time $t$  and then carrying out a regular FPCA on the linear tangent space of log-mapped data. Riemannian functional principal component (RFPC) scores, eigenfunctions, and finite-truncated representations of the log-mapped data are defined on the tangent spaces and finite-truncated representations of the data on the original manifold 
are then obtained by applying exponential maps to the  log-mapped finite-truncated data. We develop implementation and theory for  RFPCA and provide additional discussion for the  important special case where the manifold is the  Euclidean sphere, leading to Spherical Principal Component Analysis (SFPCA), in \autoref{s:model} below, where also estimation methods are introduced. The proposed SFPCA differs from  existing methods of principal component analysis on spheres \citep[e.g.,][]{jung:12,huck:16}, as these are  not targeting  functional data that consist of a sample of  time-dependent trajectories.

Theoretical properties of the proposed RFPCA are discussed in \autoref{s:theory}. \autoref{prop:resVar} states that the residual variance for a certain finite-dimensional time-varying geodesic manifold representation under the geodesic distance is upper bounded by the $L^2$ residual variance of the log-mapped data. The classical $L^2$ residual variance can be easily calculated and provides a convenient  upper bound of the residual variance under the geodesic distance. A  uniform central limit theorem for Riemannian manifold-valued functional data is presented in \autoref{thm:muCLT}. \autoref{cor:muRt} and \autoref{thm:covRt} provide asymptotic supremum convergence rates of the sample-based estimates of the  mean function, covariance function, and eigenfunctions to their population targets under proper metrics, and  the convergence rate for the sample FPC scores to their population targets is in 
\autoref{thm:xiRt}.  We also provide a consistency result for selecting the number of components used according to a criterion that is analogous to the  fraction of variance explained (FVE) criterion in  \autoref{cor:FVE}. All proofs are in the  Appendix. 

An important application for SFPCA is the  principal component analysis for longitudinal compositional data, which we will  introduce in \autoref{s:comp}, where we show that longitudinal compositional data can be mapped to functional trajectories that lie on a Euclidean sphere. We demonstrate a specific application for longitudinal compositional data in \autoref{s:data} for  behavioral patterns for fruit flies that are mapped to $S^4$, where we show that the proposed SFPCA outperforms conventional FPCA.  A second example concerns a sample of flight trajectories from Hong Kong to London, which are functional data on $S^2$. In this second example SFPCA also outperforms more conventional approaches and illustrates  the interpretability of the proposed RFPCA.  For the flight trajectory example, we demonstrate that the FPC scores produced by the  RFPCA encode more information  for classification purposes than those obtained by the classical FPCA in an $L^2$ functional space. These data examples are complemented by simulation studies reported in \autoref{s:sim}.

\section{Functional principal component analysis for random trajectories on a Riemannian manifold} \label{s:model}
\subsection{Preliminaries} \label{ss:prelim}
We briefly review the basics of Riemannian geometry essential for the study of  Riemannian manifold-valued functions; for further details, see, e.g., 
\cite{chav:06}. 
For  a smooth manifold $\cM$ with dimension $d$ and tangent spaces $T_p\cM$  at $\pinM$, a  \emph{Riemannian metric} on $\cM$ is a family of inner products $g_p: T_p\cM \times T_p\cM \rightarrow \bbR$ that varies smoothly over $\pinM$. Endowed with this Riemannian metric, $(\cM, g)$ is a \emph{Riemannian manifold}. The \emph{geodesic distance} $\dM$ is the metric on $\cM$ induced by $g$. A \emph{geodesic} is a locally length minimizing curve. The \emph{exponential map}  at $\pinM$ is defined as $\exp_p(v) = \gamma_v(1)$ where $v \in T_p\cM$ is a tangent vector at $p$, and $\gamma_v$ is a unique geodesic with initial location $\gamma_v(0) = p$ and velocity $\gamma_v'(0) = v$. If $(\cM, \dM)$ is a complete metric space, then $\exp_p$ is defined on the entire tangent space $T_p\cM$. 
The exponential map $\exp_p$ is a diffeomorphism in a neighborhood of the origin of the tangent space; the \emph{logarithm map} $\log_p$ is the inverse of $\exp_p$. 
The \emph{radius of injectivity} $\inj_p$ at $\pinM$ is the radius of the largest ball about the origin of $T_p\cM$, on which $\exp_p$ is a diffeomorphism (\autoref{fig:demo}, left panel).
If $\cN$ is a submanifold of $\cM$ with Riemannian metric $h_p:T_p\cN\times T_p\cN \rightarrow \bbR$, $(u,v) \mapsto g_p(u,v)$ for $u,v\in T_p\cN$ induced by $g$, then $(\cN, h)$ is a \emph{Riemannian submanifold} of $(\cM, g)$.

We consider a $d$-dimensional complete Riemannian submanifold $\cM$ of a Euclidean space $\bbRdz$ for $d \le d_0$, with a geodesic distance $\dM$ on $\cM$ induced by the Euclidean metric in $\bbRdz$, and a probability space  $(\Omega, \cA, P)$   with sample space $\Omega$, $\sigma$-algebra $\cA$, and probability measure $P$. With $\cX =\{x: \cT \rightarrow \cM \mid x \in \cC(\cT) \}$ denoting  the sample space of all $\cM$-valued continuous functions on a compact interval $\cT \subset \bbR$ and  $\cB(\mathcal{V})$   the Borel $\sigma$-algebra of a space $\mathcal{V}$,  the  $\cM$-valued random functions $X(t, \omega)$ are  $X: \cT \times \Omega \rightarrow \cM$, such that $X(\cdot, \omega) \in \cX$. Here  $\omega \mapsto X(\cdot, \omega)$ and $X(t, \cdot)$ are measurable with respect to $\cB(\cX)$ and $\cB(\cM)$, respectively, with $\cB(\cX)$
generated by the supremum metric $\dX: \cX \times \cX \rightarrow \bbR$, $\dX(x, y) = \supt \dM(x(t), y(t))$, 
for investigating the rates of uniform convergence. 
In the following, all vectors $v$ are column vectors and we write $X(t)$, $\tinT$, 
for $\cM$-valued random functions,  
$\normR{\cdot}$ for the  Euclidean norm, and   $\bbH = \{v: \cT \rightarrow \bbRdz,\, \int_{\cT} v(t)^T v(t) dt < \infty \}$ for  the ambient $L^2$ Hilbert space of $\bbRdz$ valued square integrable functions, equipped with the  inner product $\inner{v}{u}  = \int_{\cT} v(t)^T u(t) dt$ and norm $\Vert v \Vert = \inner{v}{v}^{1/2}$ for $u,v\in \bbH$.

\subsection{Riemannian functional principal component analysis} \label{ss:RFPCA}

As intrinsic population mean function for the $\cM$-valued random function $X(t)$, we consider the intrinsic Fr\'echet mean $\muMt$ at each time point $t \in \cT$, where
\begin{equation} \label{eq:mean}
M(p, t) = E[\dM(X(t), p)^2], \quad \mu_\cM(t) = \argmin_{p \in \cM} M(p, t),
\end{equation}
and we assume the existence and the uniqueness of the Fr\'echet means $\muMt$. 
The mean function $\muM$ is continuous due to the continuity of the sample paths of $X$, as per  \autoref{prop:muCnst} below.  
One could consider an 
alternative definition for the mean function,  $\mu_G=\argmin_{\mu} F(\mu)$, where $F(\mu) = E[\int_{\cT} \dM(X(t), \mu(t))^2 dt]$, which  
coincides with $\muM$ under a continuity assumption; we work with $\muM$ in \eqref{eq:mean}, as it matches the approach in functional PCA and allows us to 
investigate uniform convergence. 
The goal of RFPCA is to  represent the variation of the infinite dimensional object $X$ around the mean function $\mu_\cM$ in a lower dimensional submanifold, in terms of a few principal modes of variation, an approach that has been successful to represent random trajectories in the Hilbert space $L^2$ \cp{cast:86, rams:05, wang:16}. 


Given an  arbitrary system of $K$ orthonormal basis functions,  $\Psi_K = \{\psi_k \in \bbH \mid \psi_k(t) \in \Tmu, \inner{\psi_k}{\psi_l} = \delta_{kl},\, k,l = 1, \dots, K\}$,  $\delta_{kl} = 1$ if $k = l$ and 0 otherwise,  with values at each time $t \in \cT$  restricted to the $d$-dimensional tangent space $\Tmu$, which we identify with  $\bbRdz$ for convenience, we define 
the $K$ dimensional time-varying geodesic submanifold 
\begin{equation} \label{eq:MK}
\cMK(\Psi_K) \coloneqq \{x \in \cX, \, x(t) = \expmu(\sumkK a_k \psi_k(t)) \text{ for } \tinT \mid a_k \in \bbR, \, k = 1, \dots, K\}.
\end{equation}
Here $\cMK(\Psi_K)$ plays an analogous role to the linear span of a set of basis functions in Hilbert space, with expansion coefficients or coordinates $a_k$. 

In the following  we suppress the dependency of $\cMK$ on the basis functions. With projections $\Pi(x, {\cM_K})$ of an $\cM$-valued function $x \in \cX$ onto  time-varying geodesic submanifolds $\cM_K$, 
\[
\Pi(x, {\cM_K}) \coloneqq \argmin_{y\in \cM_K} \intcT \dM(y(t), x(t))^2 dt,
\]
the best $K$-dimensional approximation to $X$ minimizing the geodesic projection distance  is  the geodesic submanifold that minimizes 
\begin{equation} \label{eq:truTarget}
F_S(\cMK) = E\intcT \dM(X(t), \Pi(X, \cMK)(t))^2 dt
\end{equation}
over all time-varying geodesic submanifolds generated by $K$ basis functions.

As the minimization of \eqref{eq:truTarget} is over a family of submanifolds (or basis functions), this target is difficult to implement  in practice, except for simple situations, and therefore it is expedient to target a modified version of  \eqref{eq:truTarget} by invoking  tangent space approximations. This approximation requires that  the  log-mapped random functions
\[V(t) = \logmu(X(t))\]  
are almost surely well-defined for all $t \in \cT$, which will be 
the case if trajectories  $X(t)$ are confined to stay within the radius of injectivity at $\muMt$ for all $t \in \cT$. We require this constraint to be satisfied, which will be the case for  many manifold-valued trajectory data, including  the data we present in \autoref{s:data}. 
Then $V$ is a well-defined random function that assumes its values on the linear tangent space $\Tmu$ at time $t$. Identifying $\Tmu$ with $\bbRdz$, we may regard $V$ as a random element of  $\bbH$, the $L^2$ Hilbert space of $\bbRdz$ valued square integrable functions, and thus our analysis is independent of the choice of the coordinate systems on the tangent spaces. A practically  tractable optimality criterion to obtain manifold principal components  is then to minimize
\begin{equation} \label{eq:appTarget}
F_V(\cVK) = E(\norm{V - \Pi(V,\cV_K)}^2)
\end{equation}
over all $K$-dimensional linear subspaces $\cVK(\psi_1, \dots, \psi_K) = \{\sumkK a_k \psi_k \mid a_k \in \bbR \}$ for $\psi_k \in \bbH$, $\psi_k(t) \in \Tmu$, and $k=1, \dots, K$.  Minimizing \eqref{eq:appTarget} is immediately seen to be equivalent to a multivariate functional principal component analysis (FPCA) in $\bbRdz$ \citep{chio:14}.

\begin{figure}[h!]
\includegraphics[width=0.49\textwidth]{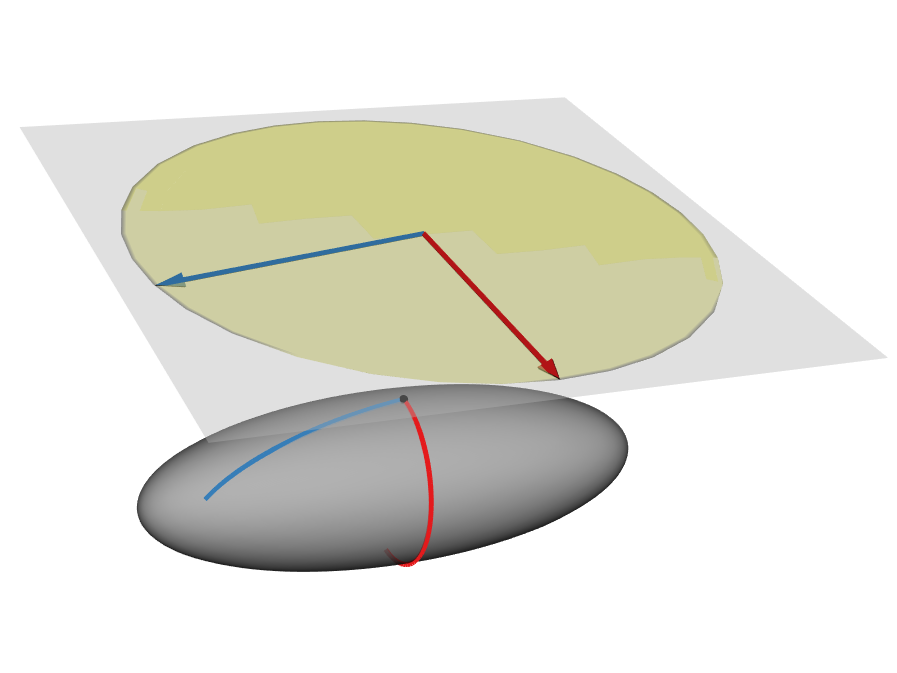}
\includegraphics[width=0.45\textwidth]{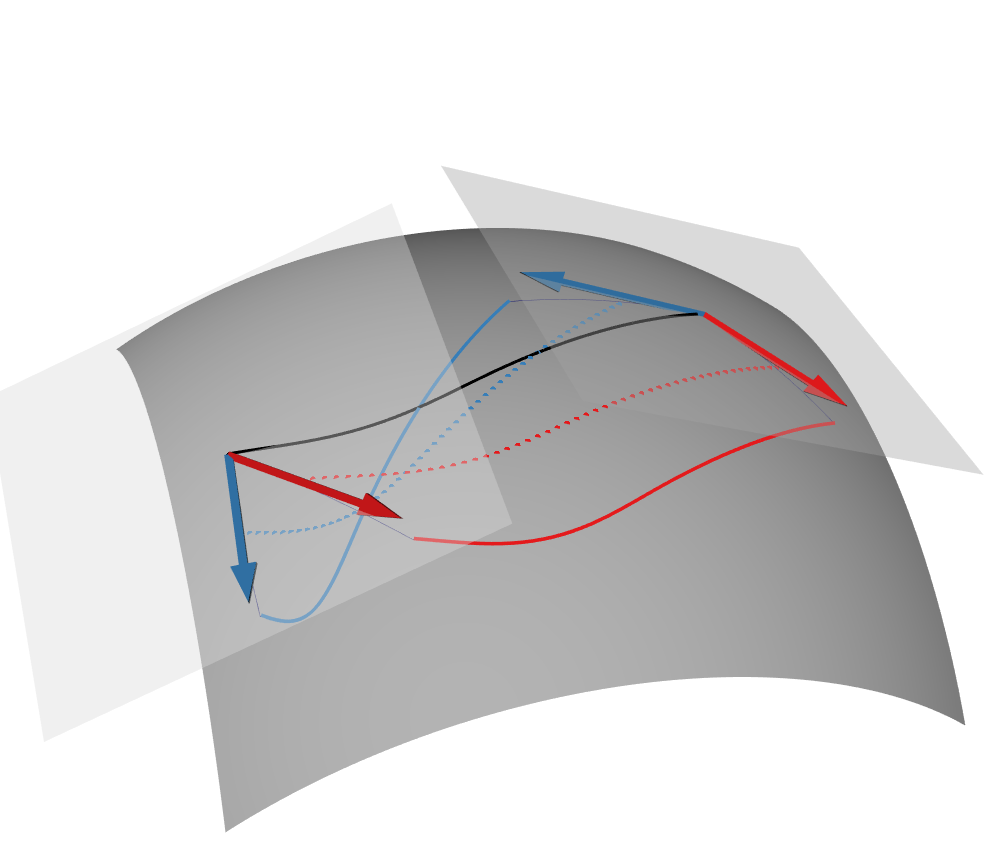}
\caption{Left panel: Two tangent vectors $v$ (red and blue arrows) in the tangent ball (yellow) centered at $p$ (black dot) with radius $\inj_p$, and their geodesics $\{\exp_p(tv)\mid t\in[0, 1]\}$ (red and blue lines). Right panel: Two trajectories $X(t)$ (red and blue solid curves), corresponding tangent vectors $V(t)$ at $t=0,1$ (arrows), and the first two eigenfunctions (red dotted, $\phi_1$, and blue dotted, $\phi_2$) mapped onto $\cM$ by the exponential maps. The red trajectory has a large score on $\phi_1$, while the blue one has a large score on $\phi_2$. The mean function is the black curve.}
\label{fig:demo}
\end{figure}

Under mild assumptions, 
the $L^2$ mean function for the log-mapped data $V(t)=\log_{\mu_\cM(t)}(X(t))$ at the Fr\'echet means is zero by Theorem~2.1 of \cite{bhat:03}. Consider the  covariance function $G$ of $V$ in the $L^2$ sense, $G:\cT\times\cT \rightarrow \bbR^{d_0^2}$, $G(t, s) = \cov(V(t), V(s)) = E(V(t)V(s)^T)$, and its associated  spectral decomposition,  $G(t, s) = \sumkinf \lambda_k \phi_k(t) \phi_k(s)^T$, 
where the $\phi_k \in \bbH: \cT \rightarrow \bbRdz$ are the orthonormal vector-valued eigenfunctions and $\lambda_k \ge 0$  the corresponding eigenvalues, for $k = 1, 2, \dots$. One obtains the  Karhunen-Lo\`eve decomposition \citep[see for example][]{hsin:15}, 
\begin{equation} \label{eq:expand}
V(t) = \sumkinf \xi_k \phi_k(t),
\end{equation}
where  $\xi_k = \int_\cT V(t) \phi_k(t) dt$ is the $k$th Riemannian functional principal component (RFPC) score, $k=1, 2, \dots$. 
A graphical demonstration of $X(t)$, $V(t)$, and $\phi_k(t)$ is in the right panel of \autoref{fig:demo}.
In practice, one can use only a finite number of components and target truncated representations of the tangent space process. Employing $K  \in \{0,1,2,\dots\}$ components,  set
\begin{equation} \label{eq:XK}
V_K(t) = \sumkK \xi_k \phi_k(t), \quad 
X_K(t) = \exp_\muMt\left(\sumkK \xi_k \phi_k(t) \right), 
\end{equation} 
where for $K=0$ the values of the sums are set to 0, so that $V_0(t) = 0$ and $X_0(t) = \mu_\cM(t)$. By classical FPCA theory, $V_K$ is the best $K$-dimensional approximation to $V$ in the sense of being the minimizing projection $\Pi(V, \cV_K)$ for  \eqref{eq:appTarget}. The truncated representation $X_K(t)$, $\tinT$ of the original $\cM$-valued random function is well-defined for $K =0, 1, \dots$ if $\cM$ is complete, by the Hopf--Rinow theorem \citep[see, e.g.,][]{chav:06}. We note that these definitions are independent of the choice of coordinate system on $\Tmu$.

To quantify how well $X_K$ approximates $X$, in analogy to  \cite{pete:16}, we define for $K = 0, 1, \dots$ the residual variance as
\begin{equation} \label{eq:U}
U_K = E \intcT \dM(X(t), X_K(t))^2 dt,
\end{equation}
and the fraction of variance explained (FVE) by the first $K$ components as
\begin{equation} \label{eq:FEV}
\FVE_K = \frac{U_0 - U_K}{U_0}. 
\end{equation}
A commonly used criterion for choosing  the number of included components $K^*$ is to select the smallest $K$ such that  FVE exceeds a specified threshold $0 < \gamma < 1$ of variance explained, 
\begin{equation}
K^* = \min\left\{K: \FVE_K \ge \gamma \right\}. \label{eq:K*}
\end{equation}
Common choices for the FVE threshold $\gamma$ are $0.9$ or $0.95$ in finite sample situations or $\gamma$  increasing with  sample size for asymptotic considerations. 

\single
\double

\subsection{Spherical functional principal component analysis} \label{ss:SFPCA}
An important special case occurs when random trajectories lie on $\cM = S^{d}$, the Euclidean sphere in $\bbR^{d_0}$ for $d_0 = d + 1$, with the Riemannian geometry induced by the Euclidean metric of the ambient space.  Then  the proposed RFPCA specializes to  spherical functional principal component analysis (SFPCA). 
We briefly review the geometry of Euclidean spheres. The geodesic distance $\dM$ on the sphere is the great-circle distance, i.e. for $p_1, p_2 \in \cM = S^{d}$
\[
\dM(p_1, p_2) = \cos^{-1}(\innerR{p_1}{p_2}).
\]

A geodesic is a segment of a great circle that connects two points on the sphere. For any point $p \in \cM$, the tangent space $T_p\cM$ is identified by $\{v \in \bbRdz \mid \innerR{v}{p} = 0\} \subset \bbRdz$, with the Euclidean inner product. Letting $\normR{\cdot}$ be the Euclidean norm in the ambient Euclidean space $\bbRdz$, then for a tangent vector $v$ on the tangent space $T_p\cM$, the exponential map is
\[
\exp_p(v) = \cos(\normR{v})p + \sin(\normR{v}) \frac{v}{\normR{v}}.
\]
The logarithm map $\log_p: \cM \setminus \{-p\} \rightarrow T_p\cM$ is the inverse of the exponential map,
\[
\log_p(q) = \frac{u}{\normR{u}} \dM(p, q),
\]
where $u = q - (\innerR{p}{q})\,p$, and $\log_p$ 
is defined everywhere with the exception of  
the antipodal point $-p$ of $p$ on $\cM$. The radius of injectivity is therefore $\pi$. The sectional curvature of a Euclidean sphere is constant.

\subsection{Estimation} \label{ss:est}
Consider a Riemannian manifold $\cM$ and $n$ independent observations $X_1, \dots, X_n$,  which are $\cM$-valued random functions that are distributed as $X$, where we assume that these functions are fully observed for $t \in \cT$. 
Population quantities for RFPCA are estimated by their empirical versions, as follows. Sample Fr\'echet means $\hmuMt$ are obtained by minimizing $M_n(\cdot, t)$ at each $\tinT$, 
\begin{equation} \label{eq:meanhat}
M_n(p, t) = \oneovern\sumin \dM(X_i(t), p)^2, \quad \hmuMt = \argmin_{p \in \cM} M_n(p, t). 
\end{equation}
We estimate the log-mapped data $V_i$ by $\hV_i(t) = \log_\hmuMt(X_i(t))$, $\tinT$; 
the covariance function $G(t,s)$ by the sample covariance function $\hG(t,s) = n^{-1} \sumin \hV_i(t) \hV_i(s)^T$ based on $\hV_i$, for $t, s \in \cT$; the $k$th eigenvalue and eigenfunction pair $(\lambda_k, \phi_k)$ of $G$ by the eigenvalue and eigenfunction $(\hlambda_k, \hphi_k)$ of $\hG$; and the $k$th RFPC score of the $i$th subject $\xi_{ik} = \int_\cT V_i(t) \phi_k(t) dt$ by $\hxi_{ik} = \int_\cT \hV_i(t)\hphi_k(t) dt$. The $K$-truncated processes  $V_{iK}$ and $X_{iK}$ for the $i$th subject $X_i$ are estimated by 
\begin{equation} \label{eq:XKhat}
\hV_{iK}(t) = \sumkK \hxi_{ik} \hphi_k(t) , \quad 
\hXiK(t) = \exp_\hmuMt \left(\sumkK \hxi_{ik} \hphi_k(t) \right),
\end{equation} 
where again for $K=0$ we set the sums to 0. 
The residual variance $U_K$ as in \eqref{eq:U}, the fraction of variance explained $\FVE_K$  as in \eqref{eq:FEV}, and the optimal $K^*$ as in \eqref{eq:K*} are respectively estimated by 
\begin{gather}
\hU_K = \oneovern \sumin \int_\cT \dM(X_i(t), \hX_{iK}(t))^2 dt, \label{eq:Uhat}\\ 
\hFVE_K = \frac{\hU_0 - \hU_K}{\hU_0}, \label{eq:FVEhat}\\
\hK^* = \min\{K:\hFVE_K \ge \gamma \}. \label{eq:hK*}
\end{gather}

Further details about  the algorithms for implementing  SFPCA can be found in the Supplementary Materials. Sometimes  functional data $X(t)$ are observed only at densely spaced time points and observations might be contaminated with measurement errors. In these situations one can presmooth the observations using smoothers that are adapted to a Riemannian manifold  \citep{jupp:87,lin:16}, treating the presmoothed curves as fully observed underlying curves. 

\section{Theoretical properties of Riemannian Functional Principal Component Analysis} \label{s:theory}


We need the following assumptions \autoref{a:rm}--\autoref{a:curvature} for the Riemannian manifold $\cM$, and \autoref{a:cont}--\autoref{a:lips} for the $\cM$-valued process $X(t)$. 
\begin{enumerate}[label=(A\arabic*)]
\item \label{a:rm} $\cM$ is a closed Riemannian submanifold of a Euclidean space $\bbRdz$, with geodesic distance $\dM$ induced by the Euclidean metric.
\item \label{a:curvature} The sectional curvature of $\cM$ is nonnegative.
\end{enumerate}
Assumption \autoref{a:rm} guarantees that the exponential map is defined on the entire tangent plane, so that $X_K(t)$ as in \eqref{eq:XK} is well-defined, while the curvature condition \autoref{a:curvature} bounds the departure between $X_K(t)$ and $X(t)$ by that of their tangent vectors. These assumptions  are satisfied for example by Euclidean spheres $S^d$. 
For the following recall $M(p,t)$ and $M_n(p,t)$ are defined as in \eqref{eq:mean} and \eqref{eq:meanhat}.

\begin{enumerate}[label=(B\arabic*)]
\item \label{a:cont} Trajectories $X(t)$ are continuous for $\tinT$ almost surely.
\item \label{a:muMExtUnq} For all $t \in \cT$, $\muMt$ and $\hmuMt$ exist and are unique, the latter almost surely.
\item \label{a:rinj} Almost surely, trajectories $X(t)$ lie in  a compact set $S_t \subset B_\cM(\muMt, r)$ for $\tinT$, where $B_\cM(\muMt,r) \subset \cM$ is an open ball centered at $\muMt$ with radius $r <  \inft\injmuMt$. 
\item \label{a:wellSep} For any $\epsilon > 0$, 
\[
\inft\infdp M(p, t) - M(\muMt, t) > 0.
\]
\item \label{a:posDef} For $v \in \Tmu\cM$, define $g_t(v) = M(\expmu(v), t)$. Then 
\[
\inf_{t \in \cT} \lambdamin(\diff{^2}{v^2}g_t(0)) > 0,
\]
where $\lambdamin(A)$ is the smallest eigenvalue of a square matrix $A$. 
\item \label{a:lips} Let $L(x)$ be the Lipschitz constant of a function $x$, i.e. $L(x) = \sup_{t \ne s} \dM(x(t), x(s)) / |t - s|$. Then $E(L(X)^2) < \infty$ and $L(\muM) < \infty$.
\end{enumerate}

Smoothness assumptions  \autoref{a:cont} and \autoref{a:lips}  for the sample paths of the observations are needed  for continuous representations, while  
existence and uniqueness of  Fr\'echet means \autoref{a:muMExtUnq} are prerequisites for an intrinsic analysis that are commonly  assumed \citep{bhat:03,pete:16:1} and depend in a complex way on the type of manifold and probability measure considered. Assumptions \autoref{a:wellSep} and \autoref{a:posDef} characterize the local behavior of the 
 criterion function $M$ around the minima and are  
standard for M-estimators \citep{bhat:17}.
Condition \autoref{a:rinj} ensures that the geodesic between $X(t)$ and $\muMt$ is unique, ensuring that the tangent vectors do not switch directions under small perturbations of the base point $\muMt$. It is satisfied for example for the sphere  $\cM = S^{d}$, if the values of the random functions are either restricted to  the positive quadrant of the sphere, as is the case for longitudinal compositional data as in \autoref{s:comp}, or if the samples are generated by $\expmu(\sumkinf \xi_k \phi_k(t))$ with bounded eigenfunctions $\phi_k$ and small scores $\xi_k$ such that $\supt|\sumkinf \xi_k \phi_k(t)| \le r$. In real data applications, \autoref{a:rinj} is justified when the $\cM$-valued samples cluster around  the intrinsic mean function, as exemplified by the flight trajectory data that  we study in  \autoref{ss:flight}. 

The following result  justifies the tangent space RFPCA approach, as the truncated representation is found to be well-defined, and the residual variance for the optimal geodesic submanifold representation bounded by that for the classical FPCA on the tangent space. 
\begin{Proposition} \label{prop:resVar}
Under \autoref{a:rm}, $X_K(t) = \expmu(V_K(t))$ is well-defined for $K = 1, 2, \dots$ and $\tinT$. If further \autoref{a:curvature} is satisfied, then
\begin{equation} 
\min_{\cM_K}\intcT \dM(X(t), \Pi(X, \cM_K)(t))^2 dt \le \intcT \dM(X(t), X_K(t))^2 dt \le \norm{V - V_K}^2.\label{eq:resVar}
\end{equation}
\end{Proposition}

The first statement  is a straightforward consequence of the Hopf-Rinow theorem, 
while the inequalities imply  that the residual variance using the best $K$-dimensional time-varying geodesic manifold approximation under geodesic distance (the left hand term) is bounded by that of the geodesic manifold produced by the proposed RFPCA (the middle term), where the latter is again bounded by the residual variance of a linear tangent space FPCA under the familiar Euclidean distance (the right hand term). The r.h.s. inequality in \eqref{eq:resVar}  affirms that  the tangent space FPCA serves as a gauge to control the preciseness of finite-dimensional approximation to the processes under the geodesic distance. An immediate consequence is that $U_K \tozero$ as $K\toinf$ for the residual variance $U_K$ in  (\ref{eq:U}), implying that  the truncated representation $X_K(t)$ is consistent for $X(t)$ when the sectional curvature of $\cM$ is nonnegative. The l.h.s. inequality gets tighter as  the samples $X(t)$ lie closer to the intrinsic mean $\muMt$, where such closeness  is not uncommon,  as demonstrated in \autoref{s:data}. The r.h.s. inequality is a consequence of the Alexandrov--Toponogov theorem for comparing geodesic triangles. 

Asymptotic properties for the estimated model components for RFPCA are studied below. 
\begin{Proposition}\label{prop:muCnst}
Under \autoref{a:rm} and \autoref{a:cont}--\autoref{a:wellSep},  $\muMt$ is continuous, $\hmuMt$ 
is continuous with probability tending to 1 as $n \toinf$, and 
\begin{equation} \label{eq:muCnst}
\supt \dM(\hmuMt, \muMt) = o_p(1).
\end{equation}
\end{Proposition}

Under additional assumptions \autoref{a:posDef} and \ref{a:lips}, the consistency in (\ref{eq:muCnst})  of the sample intrinsic mean $\hmuMt$ as an  estimator for the true intrinsic mean $\muMt$ can be strengthened through a central limit theorem on $\cC_d(\cT)$, where $\cC_d(\cT)$ is the space of $\bbR^d$-valued continuous functions on $\cT$. Let $\tau: U \rightarrow \bbR^d$ be a smooth or infinitely differentiable chart of the form $\tau(q) = \log_{p_0}(q)$, with $U=B_\cM(p_0, r_0)$, $p_0 \in \cM$, and $r_0 < \inj_{p_0}$, identifying tangent vectors in $\bbR^d$. Define chart distance $\dtau: \tau(U) \times \tau(U) \rightarrow \bbR$ by $\dtau(u, v) = \dM(\tau^{-1}(u), \tau^{-1}(v))$, its gradient $T(u,v) = [T_j(u, v)]_{j=1}^d = [\partial \dtau(u, v) /\partial v^j]_{j=1}^d$, Hessian matrix $H(u,v)$ with $(j,l)$th element $H_{jl}(u,v) = \partial^2 \dtau^2(u, v) / \partial v_j \partial v_l$, and $\Lambda(t) = E[H(\tau(X(t)), \tau(\muMt))]$. 

\begin{Theorem} \label{thm:muCLT}
Suppose that $\muMt$ and $X(t)$ are contained in the domain of $\tau$ for $t \in \cT$, the latter almost surely, and $\autoref{a:rm}$ and \autoref{a:cont}--\autoref{a:lips} hold. Then 
\begin{equation} \label{eq:muCLT}
\sqrtn [\tau(\hmuM) - \tau(\muM)] \conv{L} Z,
\end{equation}
where $Z$ is a Gaussian process with sample paths in $\cC_d(\cT)$, mean zero, and covariance $G_\mu(t, s) = \Lambda^{-1}(t) G_T(t, s) \Lambda^{-1}(s)$, where $G_T(t, s) = E[ T(\tau(X(t)), \tau(\muMt)) T(\tau(X(s)), \tau(\muMs))^T ]$, where these quantities are 
well-defined.
\end{Theorem}

\begin{Remark} \label{rmk:piece}
The first condition in \autoref{thm:muCLT} is not restrictive, since it holds at least piecewise on some finite partition of $\cT$. More precisely, due to the compactness guaranteed by \autoref{a:rm}, \autoref{a:rinj}, and  \autoref{prop:muCnst}, there exists a finite partition $\{\cT_j\}_{j=1}^N$ of $\cT$ such that $\muMt$ and $X(t)$ are contained in $B_\cM(\muM(t_j), r_j)$, for $t \in \cT_j$, $t_j \in \cM$ and $r_j < \inj_{\muM(t_j)}$, $j=1, \dots, N < \infty$. One can then define $\tau=\tau_j \coloneqq q \mapsto \log_{\muM(t_j)}(q)$ for $t \in \cT_j$ and apply \autoref{thm:muCLT} on the $j$th piece, for each $j$. 
\end{Remark}

\begin{Corollary} \label{cor:muRt}
Under $\autoref{a:rm}$ and \autoref{a:cont}--\autoref{a:lips}, 
\begin{equation} \label{eq:muRt}
\supt \dM(\hmuMt, \muMt) = O_p(n^{-1/2}).
\end{equation}
\end{Corollary}


\begin{Remark}
The intrinsic dimension $d$ is only reflected in the rate constant but not the speed of convergence. Our situation is analogous to that of  estimating the mean of Euclidean-valued random functions \citep{bosq:00}, or more generally, Fr\'echet regression with Euclidean responses \citep{pete:16:1},  where the speed of convergence does not depend on the dimension of the Euclidean space, in contrast to common nonparametric regression settings \citep{lin:16,lin:17}. The root-$n$ rate is not improvable in general since it is the optimal rate for mean estimates in the special Euclidean case. 
\end{Remark}

An immediate consequence of \autoref{cor:muRt} is the convergence of the log-mapped data. 


\begin{Corollary} \label{cor:Vi}
Under \autoref{a:rm} and \autoref{a:cont}--\autoref{a:lips}, for $i = 1, \dots, n$,
\begin{equation}
\supt \normR{\hV_i(t) - V_i(t)} = O_p(n^{-1/2}).
\end{equation}
\end{Corollary}

In the following, we use the Frobenius norm $\normF{A} = \tr(A^T A)^{1/2}$ for any real matrices $A$, 
and assume that the eigenspaces associated with positive eigenvalues $\lambda_k>0$ have multiplicity one. We obtain convergence of covariance functions, eigenvalues, and eigenfunctions on the tangent spaces, i.e., the consistency of the spectral decomposition of the sample covariance function, as follows. 

 \begin{Theorem} \label{thm:covRt}
Assume \autoref{a:rm} and \autoref{a:cont}--\autoref{a:lips} hold. Then 
\begin{gather} 
\supts \normF{\hG(t, s) - G(t, s)} = O_p(n^{-1/2}), \label{eq:covRt} \\
\sup_{k\in\bbN}|\hlambda_k - \lambda_k| = O_p(n^{-1/2}), \label{eq:lamRt}
\end{gather}
and for each $k=1, 2, \dots$ with $\lambda_k > 0$, 
\begin{gather} 
\supt \normR{\hphi_k(t) - \phi_k(t)} = O_p(n^{-1/2}). \label{eq:phiRt}
\end{gather}
\end{Theorem}

Our next result provides the  convergence rate of the RFPC scores and is a direct consequence of \autoref{cor:Vi} and \autoref{thm:covRt}.

\begin{Theorem} \label{thm:xiRt}
Under \autoref{a:rm} and \autoref{a:cont}--\autoref{a:lips}, if $\lambda_K > 0$ for some $K \ge 1$, then for each $i = 1, \dots, n$ and $k = 1,\dots, K$, 
\begin{gather} 
|\hxi_{ik} - \xi_{ik}| = O_p(n^{-1/2}), \label{eq:xiRt} \\
\supt \normR{\hV_{iK}(t) - V_{iK}(t)} = O_p(n^{-1/2}). \label{eq:VK}
\end{gather}
\end{Theorem}

To demonstrate asymptotic consistency for the number of components selected according to the FVE criterion, we consider an increasing sequence of FVE thresholds   $\gamma=\gamma_n  \uparrow 1$ as  sample size $n$ increases, which leads to a corresponding increasing sequence of $K^*=K^*_n$,
where  $K^*$ is  the smallest number of eigen-components that explains the fraction of variance  $\gamma = \gamma_n$.
One may show that  the number of components $\hK^*$ selected from the sample is consistent for the true target $K^*$ for a 
 sequence $\gamma_n$. This is formalized in the following \autoref{cor:FVE}, which is similar to Theorem~2 in \cite{pete:16}, where also specific choices of $\gamma_n$ and the corresponding sequences  $K^*$ 
were  discussed. The proof is therefore omitted. Quantities $U_0,\, U_K,\, K^*,\, \hU_0,\, \hU_K,\, \hK^*$ that appear below were  defined in \eqref{eq:U}--\eqref{eq:K*} and \eqref{eq:Uhat}--\eqref{eq:hK*}.
 
\begin{Corollary} \label{cor:FVE}
Assume \autoref{a:rm}--\autoref{a:curvature} and \autoref{a:cont}--\autoref{a:lips} hold. If the eigenvalues $\lambda_1 > \lambda_2 >\dots > 0$ are all distinct, then there exists a sequence $0 < \gamma_n \uparrow 1$ such that 
\begin{equation}
\max_{1 \le K \le K^*}\left| \frac{\hU_0 - \hU_K}{\hU_0} - \frac{U_0 - U_K}{U_0} \right| = o_p(1), \label{eq:FVEconsistent}
\end{equation}
and therefore
\begin{equation}
P(\hK^* \ne K^*) = o(1).
\end{equation}
\end{Corollary}

\section{Longitudinal compositional data analysis} \label{s:comp}
Compositional data represent proportions and are characterized by a vector $\by$ in the simplex 
\[
\cC^{J-1} = \{ \by = [y_1, \dots, y_J] \in \bbR^J \mid y_j \ge 0,\,  j=1, \dots, J;\, \sum_{j=1}^J y_j = 1\},
\]
requiring that the nonnegative proportions of all $J$ categories sum up to one. Typical examples include the geochemical composition of rocks or other data that consist of recorded percentages.  Longitudinal compositional data arise when the compositional data for the same subject are collected repeatedly at different time points. If compositions are monitored  continuously,  each sample path of longitudinal compositional data is a function $y: \cT \rightarrow \cC^{J-1}$. Analyses of such data, for example from a prospective ophthalmology study  \citep{qiu:08} or the surveillance of the composition of antimicrobial use over time \citep{adri:11}, have drawn both methodological and practical interest, but as of yet there exists no unifying methodology for  longitudinal compositional data, to the knowledge of the authors.

A direct application of standard Euclidean space methods, viewing  longitudinal compositional data as unconstrained functional data vectors \cp{chio:14},   would ignore  the non-negativity and unit sum constraints and therefore the resulting multivariate FPCA representation moves outside of the space of compositional data, diminishing the utility  of such simplistic approaches. There are various transformation that have been proposed over the years for the  analysis of compositional data to enforce the constraints, for example log-ratio transformations such as $\log(y_j/y_J)$ for $j = 1, \dots, J-1$, after which the data are  treated as Euclidean data \cp{aitc:86}, which induces the Aitchison geometry on the interior of the simplex $\cC^{J-1}$. However, these transformations cannot be defined when some of the elements in the composition are zeros, either due to the discrete and noisy nature of the observations or when the true proportions do contain actual zeros, as is the case in the  fruit fly behavior pattern data that we study in  \autoref{ss:fly} below. 

We propose to view longitudinal compositional data  as a special case of  multivariate functional data under constraints, specifically as realizations of a compositional process over time, 
\begin{equation} 
Y(t) \in \{[Y_1(t), \dots, Y_J(t)] \in \bbR^J \mid Y_j \in L^2(\cT),\, Y_j(t) \ge 0,\,  j=1, \dots, J;\, \sum_{j=1}^J Y_j(t) = 1\}, \label{eq:comp}
\end{equation}
where the component functions will also be assumed to be continuous on their domain $\cT$. 
To include the entire simplex $\cC^{J-1}$ in our longitudinal compositional data analysis, we apply square root transformations to the longitudinal compositional data $Y(t) = [Y_1(t), \dots, Y_J(t)]$, obtaining 
\begin{equation} 
X(t) = [X_1(t), \dots, X_J(t)] = [Y_1(t)^{1/2}, \dots, Y_J(t)^{1/2}]. \label{eq:scomp}
\end{equation}
A key observation is that the values of  $X(t)$ lie  on the positive quadrant of a hypersphere $S^{J-1}$ for $\tinT$, as  $X_j(t) \ge 0$ and $\sumjJ X_j(t)^2 = 1$. There is no problem with zeros as with some other proposed transformations for compositional data.  It is then a natural approach to consider  a spherical geometry for the transformed data $X(t)$. A square-root transformation and the spherical geometry for non-longitudinal  compositional data were previously considered by \cite{huck:16}. Now, since $X(t)$ assumes its values  on a quadrant of the sphere  $S^{J-1}$,  processes $X(t)$ fall into the framework of the proposed SFPCA, as described in  \autoref{ss:SFPCA}. 

Concerning the theoretical properties of SFPCA of longitudinal compositional data, the conditions on the Riemannian manifold $\cM$ needed for RFPCA are easily seen to be satisfied, due to the geometry of the Euclidean sphere and the positive quadrant constraint. We conclude 
\begin{Corollary} \label{cor:comp}
Under \autoref{a:cont} and \autoref{a:wellSep}--\autoref{a:lips}, Propositions \ref{prop:resVar} and \ref{prop:muCnst}, Theorems \ref{thm:muCLT}--\ref{thm:xiRt}, and Corollaries \ref{cor:muRt}--\ref{cor:FVE} hold for the Spherical Functional Principal Component Analysis (SFPCA)  of longitudinal compositional data $X(t)$ in  (\ref{eq:scomp}). 
\end{Corollary}

\section{Data applications} \label{s:data}
\subsection{Fruit fly behaviors} \label{ss:fly}
To illustrate the proposed SFPCA based  longitudinal compositional data analysis, we consider the lifetime behavior pattern data of \emph{D. melanogaster} \cite[common fruit fly,][]{mull:06:1}. The behavioral patterns of each fruit fly was observed instantaneously 12 times each day during its entire lifetime, and for each observation one of the five behavioral patterns, feeding, flying, resting, walking, and preening, was recorded. We analyzed the behavioral patterns in the first 30 days since eclosion for $n=106$ fruit flies with uncensored observations, aiming  to characterize and represent age-specific behavioral patterns of individual fruit flies. 
For each fruit fly, we observed the behavioral counts $[Z_1(t), \dots, Z_5(t)]$ for the five behaviors at time $t\in \cT = [0, 30]$, where the time unit is day since eclosion, and $\sum_{j=1}^5 Z_j(t)=12$ is the constrained total number of counts at each time $t$, with $0 \le Z_j(t) \le 12$ for each $j$ and $t$.  Since the day-to-day behavioral data are  noisy, we presmoothed the counts  $Z_j(t)$ of the $j$th behavior pattern over time for $j=1, \dots, 5$,  using a Nadaraya--Watson kernel smoother \citep{nada:64,wats:64} with Epanechnikov kernel and a bandwidth of five days. The  smoothed data were subsequently divided by the sum of the smoothed component values at each $t$, yielding a functional vector  $Y(t)= [Y_1(t), \dots, Y_5(t)],$ with  $Y_j(t)\ge 0$ for all $j$ and $t$ and $\sum_{j=1}^5 Y_j(t)=1$ for $t \in \cT$, thus corresponding to longitudinal compositional data.  

\begin{figure}[h!]
\centering
\includegraphics[width=.8\textwidth]{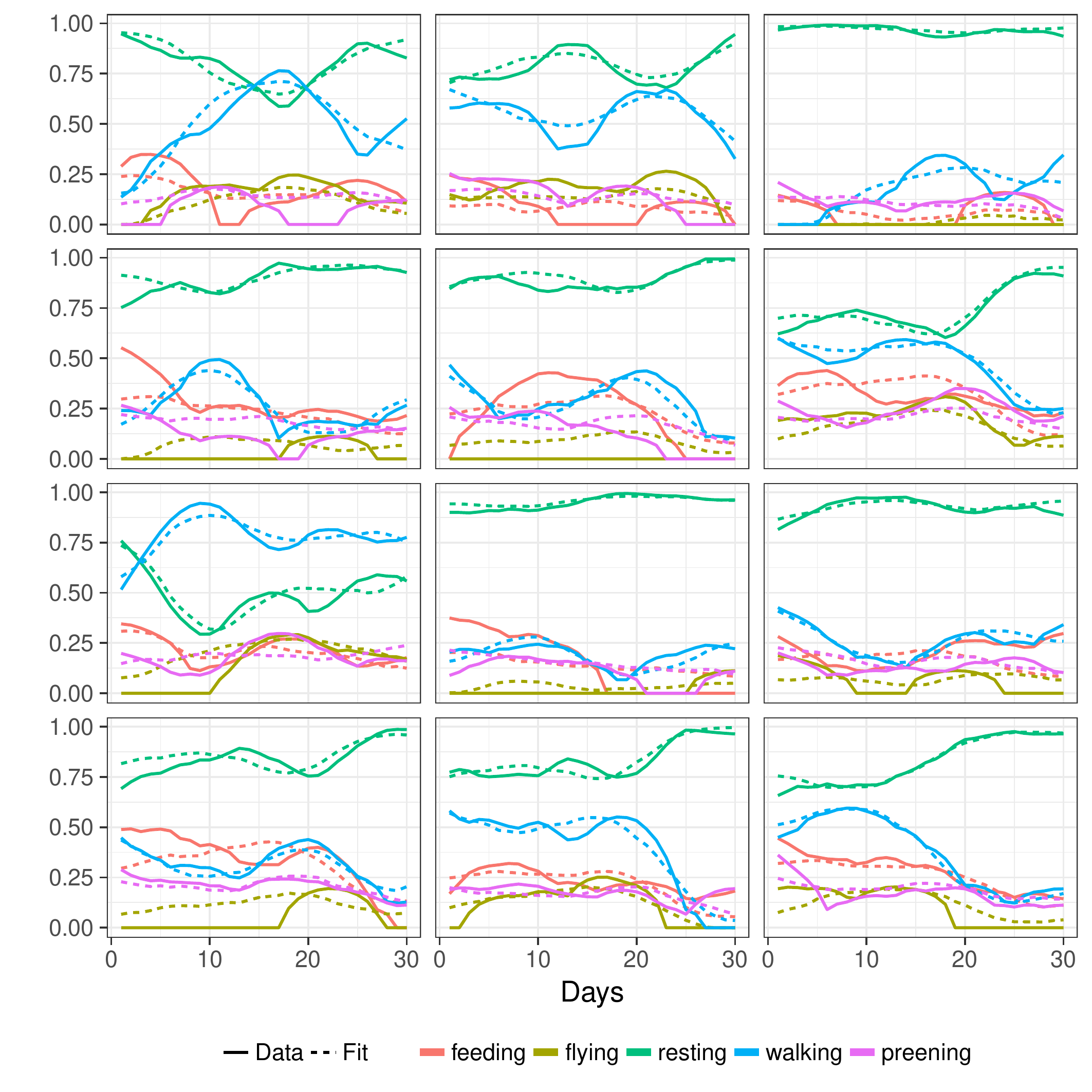}
\caption{The original data (solid lines) and SFPCA fitted trajectories (dashed lines) for 12 randomly selected fruit flies, for $K=5$ selected components.} 
\label{fig:flyFit}
\end{figure}

\begin{figure}[h!]
\includegraphics[width=.75\textwidth]{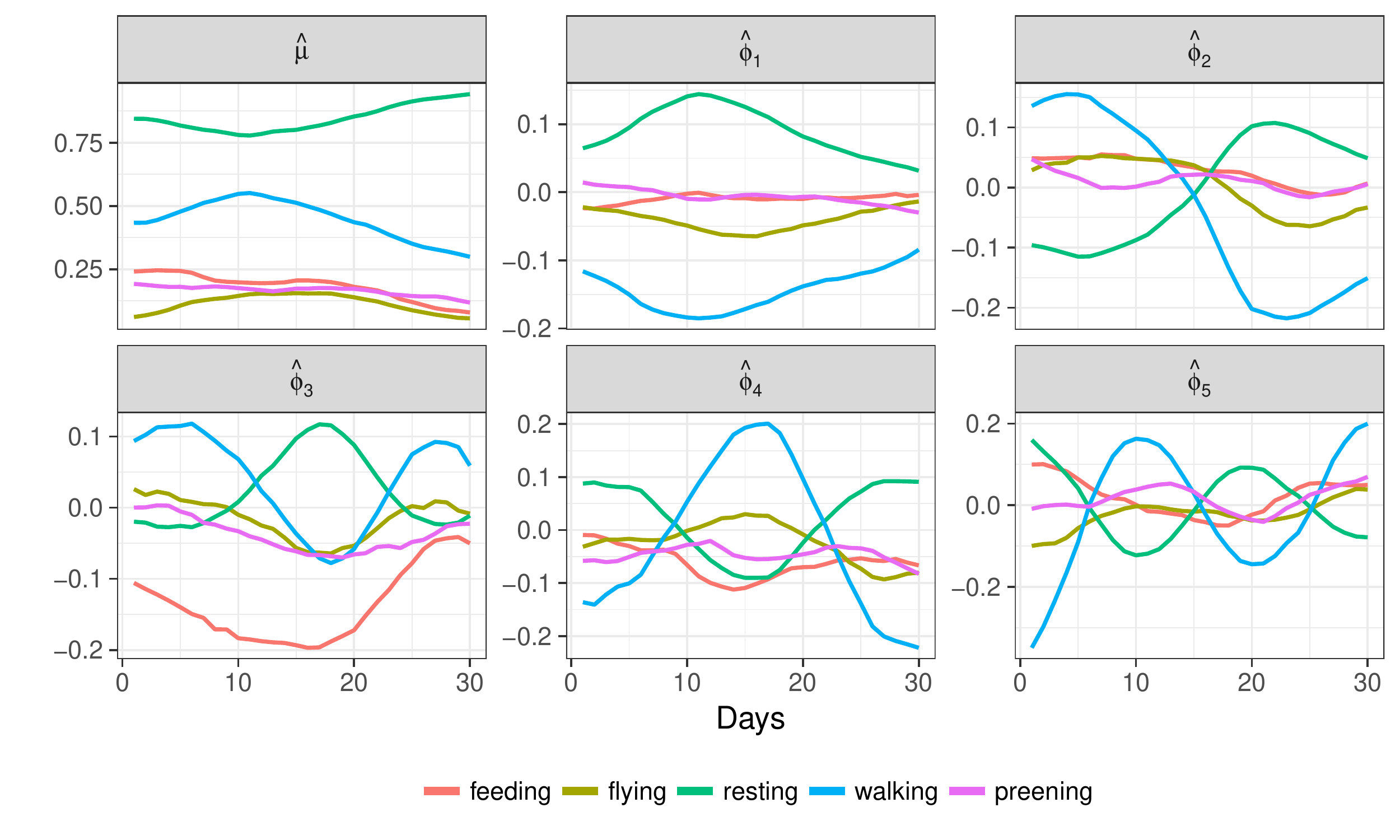}
\caption{The estimated mean functions $\hat{\mu}$  and the first five estimated spherical eigenfunctions $\hat{\phi}_1$ to $\hat{\phi}_5$ for the fly data, which together explain $81.7\%$ of the total variation. The components explain, respectively, 51.7\%,  15.0\%, 6.5\%, 5.2\%, and $3.4\%$.}
\label{fig:flyPhi}
\end{figure}

Following the approach described in \autoref{s:comp}, we model the square-root composition proportions 
$X(t) = [Y_1(t)^{1/2}, \dots, Y_5(t)^{1/2}]$ 
with  SFPCA.  The trajectories $X(t)$ and the fitted trajectories for 12 randomly selected fruit flies by SFPCA with $K=5$ components are demonstrated in \autoref{fig:flyFit}, and the mean function and the first five eigenfunctions of the corresponding SFPCA  in \autoref{fig:flyPhi}. While  resting and walking behaviors were often observed,  flying and preening occurred more rarely. SFPCA with $K = 5$ components explains 81.7\% of total variation and is seen to provide a reasonable fit to the data. 
The eigenfunctions obtained from SFPCA have a  natural interpretation: The first eigenfunction corresponds to the overall contrast of resting and moving (mainly flying and walking) over all days of observation; the second eigenfunction is a contrast of all behaviors in the early (0--15 days) and the late (16--30 days) periods; and the third eigenfunction mainly reflects  the intensity of the feeding behavior in the first 25 days.

The fraction of variance explained by the first $K$ components (FVE)  as in \eqref{eq:FVEhat} for  SFPCA and for  $L^2$~FPCA is  in \autoref{tab:flyFVE}, where  $L^2$~FPCA is  conventional multivariate FPCA \citep{rams:05}, which ignores the compositional constraints. 
The proposed SPFCA has larger FVE given any number of included components $K$. It is seen to be more parsimonious than $L^2$~FPCA and it respects the compositional constraints, in contrast to conventional FPCA. 
To explain 95\% of total variation, 14 components are needed for SFPCA, but 18 for $L^2$~FPCA.

\single
\double

\single 
\double

\single
\begin{table}[h!]
\centering
\caption{FVE (\%) by the first $K$ components for the fruit fly data.} 
\label{tab:flyFVE}
\begin{tabular}{c|ccccccccc}
$K$ & 1 & 2 & 3 & 4 & 5 & 10 & 15 & 20 & 25 \\ 
  \hline
SFPCA & 51.7 & 66.7 & 73.1 & 78.3 & 81.7 & 91.8 & 96.4 & 98.4 & 99.2 \\ 
  $L^2$ FPCA & 48.8 & 62.9 & 68.3 & 71.5 & 77.3 & 87.5 & 92.7 & 96.4 & 98.0 \\ 
\end{tabular}
\end{table}
\double

\subsection{Flight trajectories} \label{ss:flight}
\begin{figure}[h!]
\centering
\includegraphics[width=.6\textwidth]{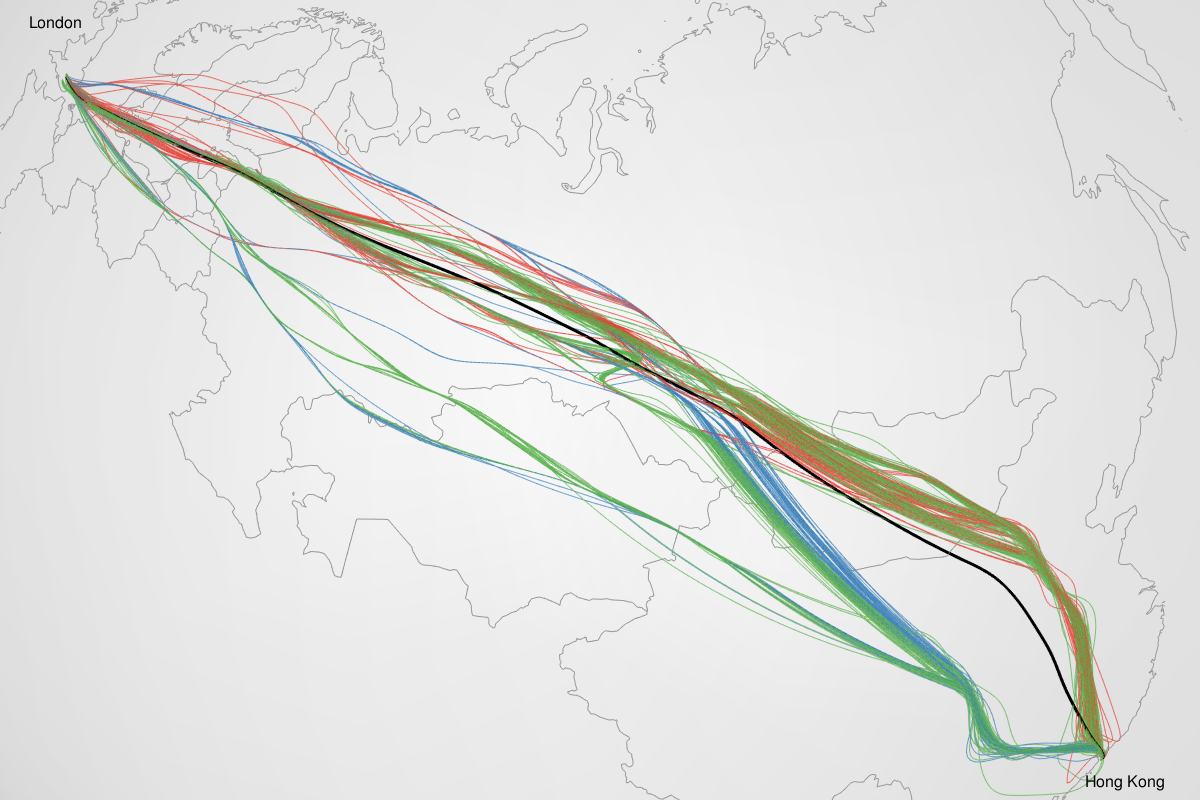}
\caption{Flight trajectories from Hong Kong to London, colored by  airline (red, British Airways; green, Cathay Pacific; blue, Virgin Atlantic), with the mean trajectory (bold black).}
\label{fig:flightSamp}
\end{figure}

A second data example concerns the trajectories of 969 commercial flights from Hong Kong to London from Jun 14, 2016 to Oct 13, 2016, of which 237 were operated  by British Airways (BAW), 612 by Cathay Pacific (CPA), and 119 by Virgin Atlantic (VIR). The data were collected from the website of FlightAware (\url{www.flightaware.com}) and  included longitude, latitude, date, and time, etc. for the whole flight, where the location was densely and accurately tracked  by ground based Automatic Dependent Surveillance--Broadcast (ADS--B) receivers. For each flight we set the takeoff time to be time 0 and the landing time to be time 1, excluding taxi time. To obtain smooth curves from the occasionally noisy data, we pre-smoothed the longitude--latitude data using kernel local linear smoothing with a very small bandwidth and then mapped the longitude--latitude trajectories onto a unit sphere $S^2$. Trajectory data of this kind on geographical spaces corresponding to the surface of the earth that may be approximated by the sphere $S^2$  have drawn extensive interest in computer science and machine learning communities  \cp{zhen:15, anir:17}. The preprocessed flight trajectories  are visualized  in \autoref{fig:flightSamp}, indicating that the flight trajectories from the three airlines overlap and are thus not easy to discriminate. We apply RFPCA in the SFPCA version to summarize and represent the flight trajectories, and to predict the operating airline based on the RFPC scores as predictors. 

\single
\double

\begin{figure}[h!]
\centering
\includegraphics[width=.6\textwidth]{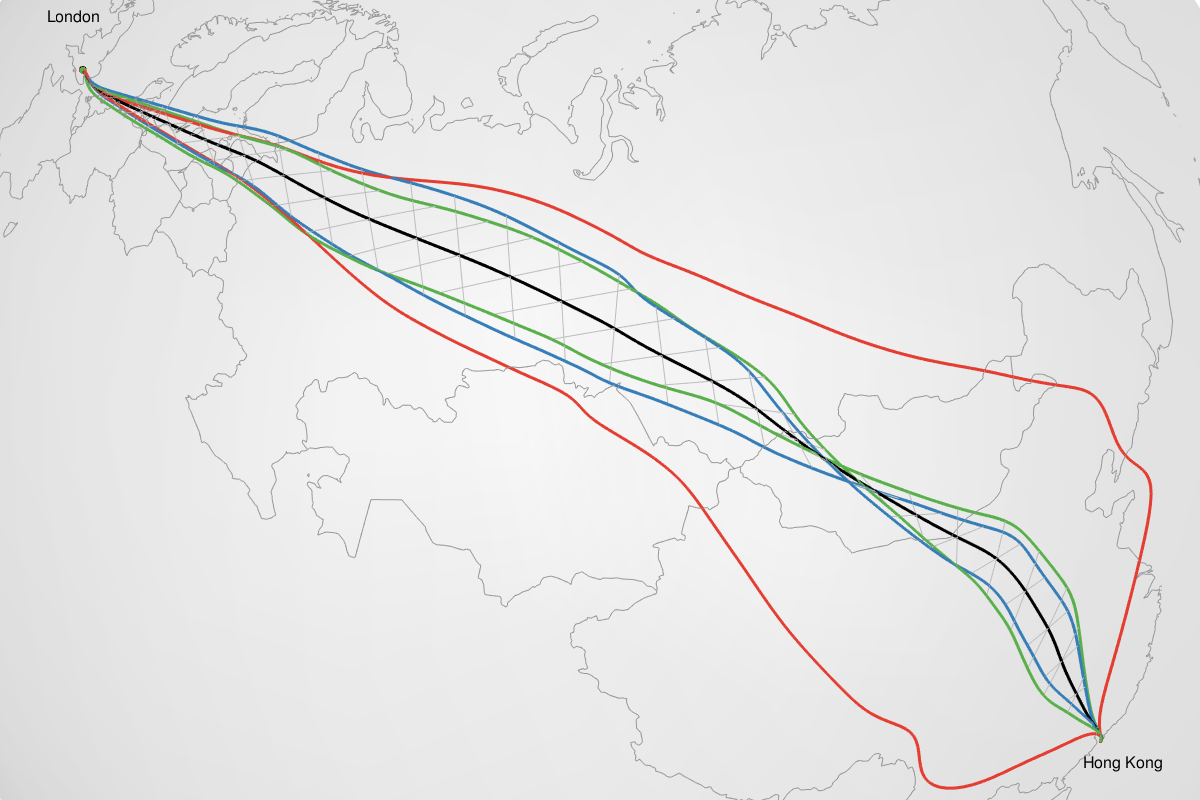}
\caption{The mean function (black) and the first three modes of variation defined as $\expmu(3\sqrt{\lambda_k}\phi_k(t))$, $k=1, 2, 3$ (red, green, and blue, respectively) produced by SFPCA. The second and the third modes of variation were joined to the time-varying mean function at a regular grid of time points to show the ``speed'' of the eigenfunctions. 
Both the second and the third eigenfunctions represent a cross from the northeast to the southwest at approximately one third of the trip, but they incorporate different speed information as shown by the thin gray lines. The first three eigenfunctions together explain in total 95\% and
each explain 72.9\%, 13.2\%, and 8.9\%, respectively,  of total variation. 
}
\label{fig:flightSpEF}
\end{figure}

The estimated mean function and the first three modes of variation obtained by SFPCA are shown in \autoref{fig:flightSpEF}, where the $k$th mode of variation is defined as $\expmu(3\sqrt{\lambda_k}\phi_k(t))$ for $k=1, 2, 3$. The first mode of variation (red) corresponds to the overall direction of deviation from the mean function (northeast vs southwest), and has roughly constant speed. We connect the second (green) and the third (blue) modes of variation and the mean function using thin gray lines at a regular grid of time in order to display speed information in the corresponding eigenfunctions. Both the second and the third eigenfunctions represent a cross from the northeast to the southwest at approximately one third of the trip, but they incorporate different speed information. The second eigenfunction encodes an overall fast trip starting to the north, while the third encodes a medium speed start to the south and then a speed up after crossing to the north. The FVE for RFPCA using the first $K=3$ eigenfunctions is 95\%, indicating a reasonably good approximation of the true trajectories. 

\single 
\double 

We next compared the FVE by the SFPCA and the \LFPCA{} for $K=1, \dots, 10$ under the geodesic distance $\dM$. Here the SFPCA was applied on the spherical data on $S^2$, while the \LFPCA{} was based on the latitude--longitude data in $\bbR^2$. A summary of the FVE for the SFPCA and the \LFPCA{} is shown in \autoref{tab:flightFVE}, using the first $K = 1, \dots, 10$ components.  Again SFPCA has higher FVE than the conventional $L^2$~FPCA for all choices of $K$, especially small $K$, where SFPCA shows somewhat  better  performance in terms of trajectory recovery. 

\single
\begin{table}[h!]
\centering
\caption{The FVE (\%) by the first $K$ components for the proposed SFPCA and the $L^2$~FPCA for the flight data.} 
\label{tab:flightFVE}
\begin{tabular}{c|ccccccccccc}
$K$ & 1 & 2 & 3 & 4 & 5 & 6 & 7 & 8 & 9 & 10 \\ 
  \hline
SFPCA & 72.9 & 86.1 & 95.0 & 96.3 & 97.0 & 97.7 & 98.3 & 98.7 & 99.0 & 99.2 \\ 
  $L^2$ FPCA & 71.2 & 84.9 & 94.6 & 96.1 & 96.8 & 97.4 & 98.1 & 98.4 & 98.8 & 99.1 \\ 
\end{tabular}
\end{table}
\double

We also aimed to predict the airline  (BAW, CPA, and VIR) from an observed flight path  by  
feeding the  FPC scores  obtained from either the proposed SFPCA or from the traditional  $L^2$~FPCA
into different multivariate classifiers, including   linear discriminant analysis (LDA),  logistic regression, and support vector machine (SVM) with radial basis kernel. For each of 200 Monte Carlo runs, we randomly selected 500 flights as training set for training and tuning  and used the rest as  test set to evaluate classification performance.  The number of components $K$ for each classifier was either fixed at  10, 15, 20, 25, 30, or selected by five-fold cross-validation (CV). The results for prediction accuracy are  in \autoref{tab:flightClass}. The SFPCA based classifiers performed better or at least equally well as the \LFPCA{} based classifiers for nearly all choices of $K$ and classifier,  where among the classifiers  
SVM performed best. 

\single 
\begin{table}[h!]
\centering
\caption{A comparison of airline classification accuracy (\%) from observed flight trajectories, using the first $K$ components for SFPCA and $L^2$ FPCA (columns),   
with $K$ either fixed or chosen by CV, for various  classifiers (rows).  All standard errors for the accuracies are below 0.12\%. The numbers in parenthesis are the number of components chosen by CV. S stands for  SFPCA and L for $L^2$~FPCA; LDA, linear discriminant analysis; MN, multinomial logistic regression; SVM, support vector machine.} 
\label{tab:flightClass}
\small
\begin{tabular}{c|cc|cc|cc|cc|cc|cccc}
& \multicolumn{2}{|c|}{$K=10$} & \multicolumn{2}{c|}{$K=15$} & \multicolumn{2}{c|}{$K=20$} & \multicolumn{2}{c|}{$K=25$} & \multicolumn{2}{c|}{$K=30$} & \multicolumn{2}{c}{$K$ chosen by CV} \\
& S & L & S & L & S & L & S &  L & S & L & S & L \\

  \hline
LDA & 76.9 & 75.8 & 79.6 & 78.4 & 81.9 & 81.5 & 82.7 & 82.5 & 83.5 & 82.3 & 83.2 (28.0) & 82.2 (26.2) \\ 
  MN & 78.5 & 76.0 & 81.8 & 79.4 & 83.8 & 82.7 & 84.6 & 84.0 & 85.2 & 83.6 & 84.8 (27.5) & 83.7 (25.7) \\ 
  SVM & 82.3 & 80.9 & 84.3 & 82.5 & 86.3 & 85.2 & 86.1 & 86.2 & 86.3 & 85.7 & 86.2 (24.6) & 85.8 (25.0) \\ 
\end{tabular}
\end{table}
\double

\section{Simulations} \label{s:sim}
To investigate the performance of trajectory recovery for the proposed RFPCA, we considered two scenarios of Riemannian manifolds: The Euclidean sphere $\cM=S^2$ in $\bbR^3$, and the special orthogonal group  $\cM = $ SO(3) of $3\times3$ rotation matrices, viewed as a Riemannian submanifold of $\bbR^{3\times3}$. We compared three approaches: The Direct (D) method, which directly optimizes \eqref{eq:truTarget} over all time-varying geodesic submanifolds $\cMK$ and therefore serves as a gold standard, implemented through discretization; the proposed RFPCA method (R) and the classical $L^2$ FPCA method (L), which ignores the Riemannian geometry. In the direct method, the sample curves and time-varying geodesic submanifolds are discretized onto a grid of 20 equally-spaced time points, and a quasi-Newton algorithm is used to maximize the criterion function \eqref{eq:truTarget}. We used FVE as our evaluation criterion, where models were fitted using $n = 50$ or 100 independent samples. 

We briefly review the Riemannian geometry for the special orthogonal group $\cM = $ \SON. The elements of $\cM$ are $N\times N$ orthogonal matrices with determinant 1, and the tangent space $T_p\cM$ is identified with the collection of $N \times N$ skew-symmetric matrices. For $p, q \in \cM$ and skew-symmetric matrices $u, v \in T_p\cM$, the Riemannian metric is $\inner{u}{v} = \tr(u^T v)$ where $\tr(\cdot)$ is the matrix trace; the Riemannian exponential map is $\exp_p(v) = \Exp(v)p$ and the logarithm map is $\log_p(q) = \Log(qp^{-1})$, where $\Exp$ and $\Log$ denote the matrix exponential and logarithm; the geodesic distance is $\dM(p,q) = \normF{\Log(qp^{-1})}$. For $N = 3$, the tangent space $T_p\cM$ is 3-dimensional and can be identified with $\bbR^3$ through \cp{chav:06} 
$
\iota: \bbR^3 \rightarrow T_p\cM,$ 
$\iota(a, b, c) = [0, -a, -b;\, a, 0, -c;\, b, c, 0]. 
$

\begin{figure}[h!]
\centering
\includegraphics[width=.29\textwidth]{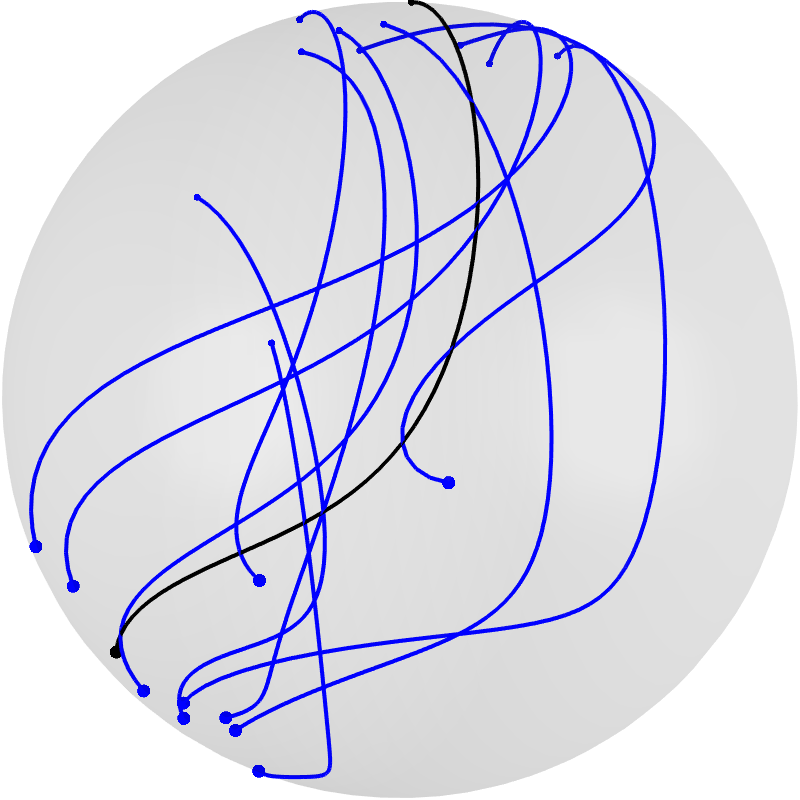}
\includegraphics[width=.29\textwidth]{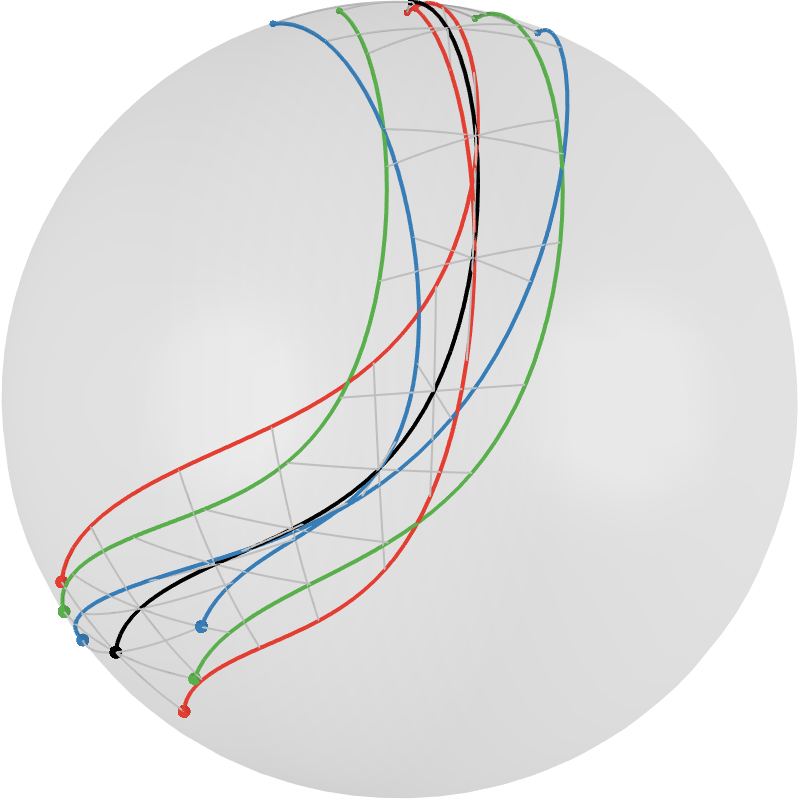}
\caption{Left: Ten randomly generated samples (dark blue) for $\cM=S^2$. Right: The first three eigenfunctions (red, green, and blue, respectively) multiplied by 0.2 and then exponentially mapped from the mean function (solid black). Light gray lines connect the mean function and the eigenfunctions at 10 equally spaced time points. Small dots denote $t=0$ and large dots  $t=1$.}
\label{fig:simSamp}
\end{figure}

The sample curves $X$ were generated as $X: \cT=[0, 1] \rightarrow \cM$, $X(t) = \expmu(\sum_{k=1}^{20}\xi_k \phi_k(t))$, with mean function $\muM(t) = \exp_{[0,0,1]}(2t, 0.3\pi\sin(\pi t), 0)$ for $\cM = S^2$, and  $\muM(t) = \Exp(\iota(2t, 0.3\pi\sin(\pi t), 0))$ for $\cM = $ SO(3). For $k = 1, \dots, 20$, the RFPC scores $\xi_k$ were generated by independent Gaussian distributions with mean zero and variance $0.07^{k/2}$. The eigenfunctions were $\phi_k(t) = 2^{-1/2} R_t [\zeta_k(t / 2), \zeta_k((t + 1)/2), 0]^T$ for $\cM = S^2$ and $\phi_k(t) = 6^{-1/2}\iota(\zeta_k(t / 3), \zeta_k((t + 1)/ 3), \zeta_k((t+2)/3))$ for $\cM = $ SO(3), $t \in [0,1]$, where $R_t$ is the rotation matrix from $[0, 0, 1]$ to $\muMt$, and $\{\zeta_k\}_{k=1}^{20}$ is the orthonormal Legendre polynomial basis on $[0,1]$. A demonstration of ten sample curves, the mean function, and the first three eigenfunctions for $\cM=S^2$ is shown in \autoref{fig:simSamp}. 

We report the mean FVE by the first $K=1, \dots, 4$ components for the investigated FPCA methods in \autoref{tab:sim}, as well as the running time, based on 200 Monte Carlo repeats. The true FVEs for $K=1, \dots, 4$ components 
were 73.5\%, 93.0\%, 98.1\%, and 99.5\%, respectively. The proposed RFPCA method had higher FVE and thus outperformed the $L^2$ FPCA in all scenarios and for all $K$, which is expected since RFPCA takes into account the curved geometry. This advantage leads to a more parsimonious representation, e.g., in the $\cM = S^2$ and $n=100$ scenario, the average $K$ required by RFPCA to achieve at least FVE$>0.95$ is one less than that for $L^2$~FPCA. The performance advantage of RFPCA over \LFPCA{} is larger for $\cM=S^2$ than for $\cM= $ SO(3),  since the former has larger sectional curvature (1 vs $1/8$). The Direct method was as expected  better than RFPCA (also for SO(3), which is not explicit in the table due to rounding), since the former optimizes the residual variation under the geodesic distance, the true target, while the latter uses the more tractable surrogate residual variation target \eqref{eq:appTarget} for $L^2$ distance on the tangent spaces. 

Each experiment was run using a single processor (Intel Xeon E5-2670 CPU @ 2.60GHz) to facilitate comparisons. Both  RFPCA and  $L^2$ FPCA are quite fast in the and take only a few seconds, though RFPCA is 1.5--3 times slower, depending on the Riemannian manifold $\cM$. The Direct method, however, was several magnitudes slower than RFPCA, due to the unstructured optimization problem, while for RFPCA spectral decomposition provides an effective solution.  The slim performance gain for the Direct method as compared to RFPCA does not justify the huge computational effort. 

\single 
\double

\single 
\begin{table}[ht]
\centering
\caption{A comparison of mean FVE (\%) and running time in the simulation study. D, direct optimization of \eqref{eq:truTarget} through discretization; R, RFPCA; L,  $L^2$~FPCA. The standard errors of the FVEs for all three methods 
were below 0.32\%.}
\label{tab:sim}
\small
\begin{tabular}{cc|ccc|ccc|ccc|ccc|ccc}
	      &     & \multicolumn{3}{|c|}{$K=1$} & \multicolumn{3}{|c|}{$K=2$} & \multicolumn{3}{|c|}{$K=3$} & \multicolumn{3}{|c}{$K=4$} & \multicolumn{3}{|c}{Time (seconds)} \\
	$\cM$ & $n$ &  D  &  R  &     L     &  D  &  R  &     L     &  D  &  R  &     L     &  D  &  R  &    L     & D  &  R  &         L          \\ \hline
	$S^2$ & 50  & 74.3 & 74.1 &     71.4      & 93.0 & 92.9 &     89.6      & 98.1 & 97.9 &     93.8      & 99.5 & 99.2 &     97.5     & 5e3 & 0.72 &          0.24          \\
	      & 100 & 74.0 & 73.8 &     70.9      & 92.9 & 92.8 &     89.2      & 98.0 & 97.9 &     93.1      & 99.4 & 99.2 &     97.3     & 1e4 & 1.01 &          0.38          \\ \hline
	SO(3) & 50  & 73.1 & 73.1 &     72.2      & 92.8 & 92.8 &     91.6      & 98.1 & 98.1 &     96.3      & 99.5 & 99.5 &     98.1     & 2e3 & 3.67 &          2.46          \\
	      & 100 & 72.9 & 72.9 &     71.8      & 92.6 & 92.6 &     91.3      & 98.0 & 98.0 &     96.1      & 99.5 & 99.5 &     97.9     & 4e3 & 6.58 &          4.94
\end{tabular}
\end{table}
\double 

\section*{Appendix: Proofs} \label{s:proofs}

\begin{proof}[Proof of \autoref{prop:resVar}]
Since $\cM$ is a closed subset of $\bbRdz$ with the induced Riemannian metric by the Euclidean metric, $\cM$ is complete. By the Hopf--Rinow theorem \citep[see, e.g., ][]{chav:06}, $\cM$ is geodesically complete, i.e. for all $\pinM$, the exponential map $\exp_p$ is defined on the entire tangent space $T_p\cM$. Therefore $X_K(t) = \expmu(V_K(t))$ is well-defined.

The first inequality in \eqref{eq:resVar} holds by the definition of projection $\Pi$. The second inequality follows from Alexandrov--Toponogov theorem \citep[e.g., Theorem~IX.5.1 in][]{chav:06}, which states if two geodesic triangles $T_1$ and $T_2$ on complete Riemannian manifolds $\cM_1$ and $\cM_2$, where 
$\cM_1$ has uniformly higher sectional curvature than $\cM_2$,
 have in common the length of two sides and the angle between the two sides, then $T_1$ has a shorter third side than $T_2$. This is applied to  triangles $(X(t), \muMt, X_K(t))$ on $\cM$ and $(V(t), 0, V_K(t))$ on $\Tmu$, identified with a Euclidean space.
\end{proof}

For the following proofs we consider the set 
\begin{equation}
\cK = \overline{\bigcup_{\tinT} B_\cM(\muMt, 2r)} \subset \cM, \label{eq:cK}
\end{equation}
where $B_\cM(p, l)$ is an open $\dM$-geodesic ball of radius $l > 0$ centered at $\pinM$, and $\overline{A}$ denotes the closure of a set $A$. Under \autoref{a:cont} and \autoref{a:rinj}, $\cK$ is closed and bounded and thus is compact, with diameter $R = \sup_{p,q\in\cK} \dM(p,q)$. Then $\muMt, \hmuMt, X(t) \in \cK$ for all $\tinT$. For the asymptotic results we will consider the compact set $\cK$.

\begin{proof}[Proof of \autoref{prop:muCnst}]
To obtain the uniform consistency results of $\hmuMt$, we need to show 
\begin{gather}
\supt\suppK |M_n(p, t) - M(p, t)| = o_p(1), \label{eq:p1:MnConv} \\
\supt |M_n(\hmuMt, t) - M(\muMt, t)| = o_p(1), \label{eq:p1:MnmuMConv}
\end{gather}
and for any $\epsilon > 0$, there exist $a = a(\epsilon) > 0$ such that 
\begin{equation}
\inft \infdp [M_n(p, t) - M(\muMt, t)] \ge a - o_p(1) \label{eq:p1:MnWellSep}.
\end{equation}
Then by \eqref{eq:p1:MnmuMConv} and \eqref{eq:p1:MnWellSep}, for any $\delta > 0$, there exists $N \in \bbN$ such that $n \ge N$ implies the event
\[
E= \{\supt |M_n(\hmuMt, t) - M(\muMt, t)| \le a / 3 \} \cap \{\inft \infdp [M_n(p, t) - M(\muMt, t)] \ge 2a/3 \}
\]
holds with probability greater than $1 - \delta$. This implies that on $E$, $\supt \dM(\hmuMt, \muMt) \le \epsilon$, and therefore  the consistency of $\hmuM$. 

Proof of \eqref{eq:p1:MnConv}:  We first obtain the auxiliary result
\begin{equation} \label{eq:p1:dXtCont}
\lim_{\delta \downarrow 0} E \left[\suptsdelta \dM(X(t), X(s))\right] = 0
\end{equation}
by dominated convergence, \autoref{a:cont}, and the boundedness of $\cK$ \eqref{eq:cK}. Note that for any $p,q,w \in \cK$, 
\[
|\dM(p,w)^2 - \dM(q,w)^2| = |\dM(p,w) + \dM(q,w)| \cdot |\dM(p,w) - \dM(q,w)| \le 2R\dM(p,q)
\] 
by the triangle inequality, where $R$ is the diameter of $\cK$. Then 
\begin{align*}
\suptspq |M_n(p, t) - M_n(q, s)| & \le \suptspq |M_n(p, s) - M_n(q, s)| + \suptspq|M_n(p,t) - M_n(p, s)| \\
& \le 2R\delta + \frac{2R}{n}\sumin\suptsdelta \dM(X_i(t), X_i(s)) \\
& =  2R\delta + 2RE\left[\suptsdelta\dM(X(t), X(s))\right]  + o_p(1),
\end{align*}
where the last equality is due to the weak law of large numbers (WLLN). Due to \eqref{eq:p1:dXtCont}, the quantity in the last display can be made arbitrarily close to zero (in probability) by letting $\deltatozero$ and $\ntoinf$. Therefore, for any $\epsilon > 0$ and $\eta > 0$, there exist $\delta > 0$ such that
\[
\limsup_{\ntoinf} P(\suptspq|M_n(p,t) - M_n(q,s)| > \epsilon) < \eta,
\]
proving the asymptotic equicontinuity of $M_n$ on $\cK\times\cT$. This and the pointwise convergence of $M_n(p, t)$ to $M(p, t)$ by the WLLN imply \eqref{eq:p1:MnConv}
by Theorem 1.5.4 and Theorem 1.5.7 of \cite{vand:96}. 

Proof of \eqref{eq:p1:MnmuMConv}: Since $\hmuMt$ and $\muMt$ are the minimizers of $M_n(\cdot, t)$ and $M(\cdot, t)$, respectively, 
$
|M_n(\hmuMt,t) - M(\muMt, t)| \le \max(M_n(\muMt,t) - M(\muMt,t), M(\hmuMt,t) - M_n(\hmuMt,t)) \le \suppK|M_n(p,t)-M(p,t)|.
$
Take suprema over $\tinT$ and then apply \eqref{eq:p1:MnConv} to obtain \eqref{eq:p1:MnmuMConv}.

Proof of  \eqref{eq:p1:MnWellSep}: Fix $\epsilon > 0$ and let $a = a(\epsilon) = \inft\infdp [M(p, t) - M(\muMt, t)]> 0$. For small enough $\epsilon$, 
\begin{align*}
\inft\infdp  & [M_n(p,t) - M(\muMt, t)] = \inft\infdpK[M_n(p,t) - M(\muMt, t)]\\
& = \inft\infdpK [M(p,t) - M(\muMt, t) + M_n(p,t) - M(p,t)] \\
& \ge a - \supt\supdpK|M_n(p,t) - M(p,t)|  = a - o_p(1),
\end{align*}
where the first equality is due to $\hmuMt\in\cK$ and the continuity of $M_n$, the inequality to \autoref{a:wellSep}, and the last equality to \eqref{eq:p1:MnConv}. 
For the continuity of $\muM$, note for any $t_0, t_1 \in \cT$,
\begin{align*}
|M(\muM(t_1), t_0) - M(\muM(t_0), t_0)| & \le |M(\muM(t_1), t_1) - M(\muM(t_0), t_0)| + |M(\muM(t_1), t_0) - M(\muM(t_1), t_1)| \\
& \le \suppK|M(p,t_1) - M(p,t_0)| +  2R E[\dM(X(t_0), X(t_1))] \\
& \le 4R E[\dM(X(t_0), X(t_1))] \tozero
\end{align*}
as $t_1 \rightarrow t_0$ by \autoref{a:cont}, where the second inequality is due to the fact that $\muM(t_l)$ minimizes $M(\cdot, t_l)$ for $l=0, 1$. Then by \autoref{a:wellSep}, $\dM(\muM(t_1), \muM(t_0)) \tozero$ as $t_1 \rightarrow t_0$, proving the continuity of $\muM$. The continuity for $\hmuM$ is similarly proven by in probability arguments.
\end{proof}

\begin{proof}[Proof of \autoref{thm:muCLT}]
The proof idea is similar to that of Theorem 2.1 in \cite{bhat:05}. To lighten notations, let $Y(t) = \tau(X(t))$, $Y_i(t)=\tau(X_i(t))$, $\nu(t) = \tau(\muMt)$, and $\hnu(t) = \tau(\hmuMt)$. 
The squared distance $\dM(p, q)^2$ is smooth at $(p, q)$ if $\dM(p,q) < \inj_p$, due to the smoothness of the exponential map \citep[Theorem I.3.2]{chav:06}. Then $\dtau(u,v)^2$ is smooth on the compact set $\{(u,v) \in \tau(U) \times \tau(U) \subset \bbR^d \times \bbR^d \mid \dM(\tau^{-1}(u), \tau^{-1}(v)) \le r\}$ and thus $T(Y(t), \nu(t))$ and $H(Y(t), \nu(t))$ are well defined,  by \autoref{a:rinj} and since the domain $U$ of $\tau$ is bounded. 
Define
\begin{align}
h_t(v) = E[ \dtau(Y(t), v)^2 ], \label{eq:gt} \\
h_{nt}(v) = \frac{1}{n} \sumin \dtau(Y_i(t), v)^2. \label{eq:gnt}
\end{align}
Since $\nu(t)$ is the minimal point of \eqref{eq:gt},  
\begin{equation} 
E[T_j(Y(t), \nu(t))] = E\left[\left. \diff{}{v_j}\dtau^2(Y(t), v)  \right\rvert_{v=\nu(t)}\right] = \diff{}{v_j}h_t(\nu(t)) = 0, \label{eq:T1st}
\end{equation}
for $j=1, \dots, d$. Similarly, differentiating  \eqref{eq:gnt} and applying Taylor's theorem,
\begin{align}
0 & = \frac{1}{\sqrtn}\sumin T_j (Y_i(t), \hnu(t)) \nonumber \\
& = \frac{1}{\sqrtn} \sumin T_j(Y_i(t), \nu(t)) + \sum_{l=1}^d \sqrtn[\hnu_l(t) - \nu_l(t)] \, \oneovern \sumin H_{jl}(Y_i(t), \nu(t)) + R_{nj}(t), \label{eq:taylor}
\end{align}
where $\hnu_l(t)$ and $\nu_l(t)$ are the $l$th component of $\hnu(t)$ and $\nu(t)$, and 
\begin{equation}
R_{nj}(t) = \sum_{l=1}^d \sqrtn[\hnu_l(t) - \nu_l(t)] \oneovern\sumin \left[ H_{jl}(Y_i(t), \tnu_{jl}(t)) - H_{jl}(Y_i(t), \nu(t)) \right],
\end{equation}
for some $\tnu_{jl}(t)$ lying between $\hnu_l(t)$ and $\nu_l(t)$. 

Due to the smoothness of $\dtau^2$, \autoref{a:rinj}, and \autoref{a:lips}, for $j,l=1, \dots, d$, 
\begin{equation}
E\supt T_j(Y_i(t), \nu(t))^2  < \infty, \label{eq:T2nd} \quad 
E\supt H_{jl}(Y_i(t), \nu(t))^2 < \infty, 
\end{equation}
\begin{equation}
\lim_{\epsilon\downarrow 0 } E \supt \sup_{\norm{\theta - \nu(t)} \le \epsilon} |H_{jl}(Y(t), \theta) - H_{jl}(Y(t), \nu(t))|  = 0. \label{eq:Heps}
\end{equation}
By \autoref{a:lips}, we also have 
$
\lim_{\epsilon \downarrow 0} E \sup_{|t-s| < \epsilon} |H_{jl}(Y(t),\nu(t)) - H_{jl}(Y(s), \nu(s))| \tozero,
$
which implies the asymptotic equicontinuity of $n^{-1}\sumin H_{jl}(Y_i(t), \nu(t))$ on $\tinT$, and thus
\begin{equation}
\supt\left| \frac{1}{n}\sumin H_{jl}(Y_i(t), \nu(t)) - E[H_{jl}(Y_i(t), \nu(t))] \right| = o_p(1), \label{eq:Hunif}
\end{equation}
by Theorem 1.5.4 and Theorem 1.5.7 of \cite{vand:96}. In view of \eqref{eq:T2nd}--\eqref{eq:Hunif} and \autoref{prop:muCnst}, we may write \eqref{eq:taylor} into matrix form
\begin{equation}
[\Lambda(t) + E_n(t)] \sqrtn[\hnu(t) - \nu(t)] = - \frac{1}{\sqrtn}\sumin T(Y_i(t), \nu(t)), \label{eq:hnuCLT}
\end{equation}
where $\Lambda(t) = E[ H(Y(t), \nu(t)) ]$ and $E_n(t)$ is some random matrix with $\supt \normF{E_n(t)} = o_p(1)$. By \autoref{a:lips}, $T_j(Y_i(t), \nu(t))$ is Lipschitz in $t$ with a square integrable Lipschitz constant, so one can apply a Banach space central limit theorem \citep{jain:75}
\begin{equation}
\frac{1}{\sqrtn}\sumin T(Y_i, \nu) \conv{L} W,
\end{equation}
where $W$ is a Gaussian process with sample paths in $\cC_d(\cT)$, mean 0, and covariance $G_T(t, s) = E[T(Y(t), \nu(t))T(Y(s), \nu(s))^T]$. 

We conclude the proof by showing 
\begin{equation}
\inft \lambda_{\min}(\Lambda(t)) > 0. \label{eq:lamMinLambda}
\end{equation}
Let $\phi_t(v) = \log_{\mu_\cM(t)}(v)$, $f_t = \phi_t \circ \tau^{-1}$, and 
$g_t(v) = E[\dM(X(t), \exp_{\muMt}(v))^2]$, so $h_t(v) = g_t(f(v))$. Observe
\begin{align}
\diff{^2}{v_j\partial v_l} h_t(v) = \left( \diff{}{v_j} f_t(v) \right)^T \diff{^2}{v^2}g_t(v) \left( \diff{}{v_l} f_t(v) \right) + \diff{}{v}g_t(v)^T \diff{^2}{v_j\partial v_l} f_t(v). 
\end{align}
The second term vanishes at $v = \nu(t)$ by \eqref{eq:T1st}, so in matrix form
\begin{equation}
\Lambda(t) = \diff{^2}{v^2}h_t(\nu(t)) = \left(\diff{}{v}f_t(\nu(t)) \right)^T \diff{^2}{v^2}g_t(0) \left(\diff{}{v}f_t(\nu(t)) \right). 
\end{equation}
The gradient of $f_t$ is nonsingular at $\nu(t)$ since it is a local diffeomorphism. Then $\Lambda(t)$ is positive definite for all $\tinT$ by \autoref{a:posDef}, and \eqref{eq:lamMinLambda} follows by continuity. 

\end{proof}

\begin{proof}[Proof of \autoref{cor:muRt}]
Note $\dM(\hmuMt, \muMt) = \dtau(\hnu(t), \nu(t))$. By Taylor's theorem around $v = \nu(t)$, 
\[
\dtau(\nu(t), \hnu(t))^2 = 
[\hnu(t) - \nu(t)]^T \left[\left.\diff{}{v^2}\dtau^2(\nu(t), v)\right\rvert_{v=\tnu(t)}\right] [\hnu(t) - \nu(t)],
\]
where $\tnu(t)$ lies between $\hnu(t)$ and $\nu(t)$, since $\dtau^2(u,v)$ and $\partial\dtau^2(u, v)/\partial v$ both vanish at $u=v$. 
The result then follows from \autoref{thm:muCLT}, \autoref{rmk:piece}, and \autoref{prop:muCnst}.
\end{proof}

\begin{proof}[Proof of \autoref{cor:Vi}]
\begin{align}
\supt \normR{\hV_i(t) - V_i(t)} & = \supt \normR{\logmu(X_i(t)) - \loghmu(X_i(t))}\\
& \lesssim \supt |\dM(\hmuMt, \muMt)|,
\end{align}
where the last inequality is due to \autoref{a:rinj} and the fact that $\log_p(q)$ is continuously differentiable in $(p, q)$ \citep[Theorem I.3.2 in][]{chav:06}.
\end{proof}

\begin{proof}[Proof of \autoref{thm:covRt}]
Denote $\tG(t, s) = \oneovern \sumin V_i(t)V_i(s)^T$. Then 
\begin{align}
\supts & \normF{\hG(t, s) - G(t, s)}\le \supts\normF{\hG(t, s) - \tG(t, s)} + \supts\normF{\tG(t,s) - G(t, s)} \nonumber\\
& \le \oneovern\sumin \supts\normF{\hV_i(t)\hV_i(s)^T - V_i(t)V_i(s)^T} + \supts\normF{\oneovern\sumin V_i(t)V_i(s)^T - G(t, s)} \label{eq:p1:GTerms}
\end{align}
Since $\supts\normF{V_i(t)V_i(s)^T} < R^2$, viewing $V_i(t)V_i(s)^T$ as random elements in $L_\infty(\cT\times\cT, \bbR^{d^2})$ the second term is $O_p(n^{-1/2})$ by Theorem~2.8 in \cite{bosq:00}. For the first term, note
\begin{align*}
\normF{\hV_i(t)\hV_i(s)^T - V_i(t)V_i(s)^T} & \le \normF{(\hV_i(t) - V_i(t))\hV_i(s)^T} + \normF{V_i(t)(\hV_i(s) - V_i(s))^T} \\
& \le \normR{\hV_i(s)}\normR{\hV_i(t) - V_i(t)} + \normR{V_i(t)}\normR{\hV_i(s) - V_i(s)} \\
& \lesssim \supt \dM(\hmuMt, \muMt),
\end{align*}
where the second inequality is due to the  properties of the Frobenius norm, and the last is due to \autoref{cor:Vi} and \autoref{a:rinj}. Therefore, by \autoref{cor:muRt} the first term in \eqref{eq:p1:GTerms} is $O_p(n^{-1/2})$ and \eqref{eq:covRt} follows. Result \eqref{eq:lamRt} follows from applying Theorem~4.2.8 in \cite{hsin:15} and from the fact that the operator norm is dominated by the Hilbert-Schmidt norm. 

To prove \eqref{eq:phiRt}, Theorem~5.1.8 in \cite{hsin:15} and Bessel's inequality imply
\begin{equation} \label{eq:phiL2Rt}
\norm{\hphi_k - \phi_k} = O_p(n^{-1/2}).
\end{equation} 
Then note that for any $\tinT$, 
\begin{align*}
\normR{\hphi_k(t) - \phi_k(t)} & = \normR{ \int \frac{1}{\hlambda_k}\hG(t, s)\hphi_k(s) ds - \int \frac{1}{\lambda_k}G(t, s)\phi_k(s) ds } \\
& = \Vert(\frac{1}{\hlambda_k} - \frac{1}{\lambda_k})\int\hG(t, s)\hphi_k(s) ds + \frac{1}{\lambda_k} \int \hG(t, s)(\hphi_k(s) - \phi_k(s)) \\
& +  (\hG(t,s) - G(t,s)) \phi_k(s) ds \Vert_E \\
& = O_p\left( \left|\frac{1}{\hlambda_k} - \frac{1}{\lambda_k}\right| + \norm{\hphi_k - \phi_k} + \supts\normF{\hG(t, s) - G(t, s)}\right),
\end{align*}
which is of order $O_p(n^{-1/2})$ by \eqref{eq:lamRt}, \eqref{eq:phiL2Rt}, and \eqref{eq:covRt}. 
Since the r.h.s. does not involve $t$, taking suprema on both sides over $\tinT$ concludes the proof. 
\end{proof}

\begin{proof}[Proof of \autoref{cor:comp}]
Conditions \autoref{a:rm}--\autoref{a:curvature} and \autoref{a:rinj} hold for longitudinal compositional data analysis, while \autoref{a:muMExtUnq} holds by Theorem~2.1 in \cite{afsa:11}. 
\end{proof}


\single 
\references 

%

\newpage

\section*{Supplementary Materials}
\subsection*{S1. Algorithms for the RFPCA of Compositional Data}
The following Algorithms~\ref{algo:SFPCA}--\ref{algo:FVE} are provided as examples for RFPCA applied  to longitudinal compositional data $Z(t)$ or spherical trajectories $X(t)$. For longitudinal compositional data $Z(t)$, we initialize by defining $X(t)$ as the componentwise square root of $Z(t)$, which then lies on a Euclidean sphere $S^d$. We assume the trajectories are observed at $t = t_j = (j-1)/(m-1)$ for $j=1, \dots, m$, and all vectors are by default column vectors. Very similar algorithms for SFPCA have also been proposed by  \cite{anir:17}. 

The time complexity for \autoref{algo:SFPCA} is $O(nmf(d) + nm^2d^2 + (md)^3)$, where $f(d)$ is the cost for calculating a Fr\'echet mean in \autoref{a1:fmean}, which is typically $O(nd)$ or $O(nd^2)$ for gradient descent or quasi-Newton type optimizers per iteration, respectively. The most demanding  computational step  for the multivariate FPCA is $O(nm^2d^2)$ for \autoref{a1:G} and $O((md)^3)$ for the eigendecomposition in \autoref{a1:eig}. The computational cost for \autoref{algo:KRep} is $O(md)$ and that for \autoref{algo:FVE} is $O(nmd)$. \vspace{1cm}

\begin{algorithm}[H]
\SetKwProg{Fn}{Function}{}{end}
\SetKwFunction{RFPCA}{RFPCA}
\SetKwFunction{SFPCA}{SFPCA}
\SetKwFunction{Vec}{Vec}
\SetKwFunction{Eigen}{Eigen}%
\newcommand{\forcond}{$i=0$ \KwTo $n$}
\single
\DontPrintSemicolon
\KwData{$S^d$-valued trajectories $X_1(t), \dots, X_n(t)$}
\KwResult{$\hmuMt$, $\hV_i(t)$, $\hxi_{ik}$, $\hphi_k(t)$, $\hlambda_k$, for $i=1, \dots, n$ and $k=1, \dots, K$}
\tcp*[l]{Obtain the intrinsic mean function and tangent vectors}
\For{$j \in \{1, \dots, m\}$}{\label{a1:1start}
	$	\hmuM(t_j) \leftarrow \argmin_{p \in S^{d}} n^{-1}\sumin [\cos^{-1}(p^T X_i(t_j))]^2$ \label{a1:fmean}

	\For{$i \in \{1, \dots, n\}$}{
		$\hV_i(t_j) = \frac{u}{\sqrt{u^Tu}} \cos^{-1}(\hmuM(t_j)^T X_i(t_j))$, where $u = X_i(t_j) - (\hmuM(t_j)^TX_i(t_j))\muM(t_j)$
	}
}\label{a1:1end}

\tcp*[l]{A multivariate FPCA. Vec(A) stacks the columns of A.}
$\hbVi \leftarrow [\hV_i(t_1), \dots, \hV_i(t_m)]^T$, $\hbG \leftarrow n^{-1}\sumin \Vec(\hbVi)\Vec(\hbVi)^T$ \label{a1:G}

$[\omega,\, \bPsi] \leftarrow \Eigen(\hbVi)$, for eigenvalues $\omega=[\omega_1, \dots, \omega_m]^T$ and eigenvectors $\bPsi = [\psi_1, \dots, \psi_m]$ \label{a1:eig}

\For{$k \in \{1, \dots, K\}$}{
Write $\hbPhi_k = [\hphi_k(t_1), \dots, \hphi_k(t_m)]^T$, $\Vec(\hbPhi_k) \leftarrow m^{1/2}\psi_k$, $\hlambda_k \leftarrow m^{-1}\omega_k$, $\hxi_{ik} \leftarrow m^{-1} \Vec(\hbVi)^T\Vec(\hbPhi_k)$
}

\caption{Spherical functional principal component analysis (SFPCA)}%
\label{algo:SFPCA}
\end{algorithm}

\vspace{1cm}

\begin{algorithm}[H]
\SetKwProg{Fn}{Function}{}{end}
\SetKwFunction{KRep}{KRep}
\SetKwFunction{Vec}{Vec}
\SetKwFunction{Eigen}{Eigen}%
\newcommand{\forcond}{$i=0$ \KwTo $n$}
\single
  \caption{Truncated $K$-dimensional representations}%
  \label{algo:KRep}
\DontPrintSemicolon

\KwData{$\hmuM(t)$, $\{(\hxi_{ik}, \hphi_k(t))\}_{k=1}^K$}
\KwResult{$\hX_{iK}(t)$, $\hV_{iK}(t)$}
\For{$j \in \{1, \dots, m\}$}{
	$\hV_{iK}(t_j) \leftarrow \sumkK \hxi_{ik}\hphi_k(t_j) $
	
	$\hXiK(t_j) \leftarrow \cos(\normR{\hV_{iK}(t_j)})\hmuMt + \sin(\normR{\hV_{iK}(t_j)})\normR{\hV_{iK}(t_j)}^{-1}\hV_{iK}(t_j) $
}
\end{algorithm}

\vspace{1cm}

\begin{algorithm}
\SetKwProg{Fn}{Function}{}{end}
\SetKwFunction{KRep}{KRep}
\SetKwFunction{Vec}{Vec}
\SetKwFunction{Eigen}{Eigen}%
\newcommand{\forcond}{$i=0$ \KwTo $n$}
\single
  \caption{Calculate FVE}%
  \label{algo:FVE}
\DontPrintSemicolon

\KwData{Outputs from \autoref{algo:SFPCA}}
\KwResult{$\hFVE_K$}

$\hU_0 \leftarrow n^{-1}\sumin\int_{\cT} \dM^2(X_i(t), \hmuMt) dt$

\For{$i \in \{1, \dots, n\}$}{
	Use \autoref{algo:KRep} to obtain $\hXiK(t)$
}

$\hU_K \leftarrow n^{-1}\sumin\int_{\cT} \dM^2(X_i(t), \hXiK(t)) dt$

$\hFVE_K = (\hU_0 - \hU_K) / \hU_0$
\end{algorithm}

\subsection*{S.2 Additional simulations}
We conducted an additional simulation study to investigate the scalability of the RFPCA algorithms to higher dimensions $d$, on the unit sphere $\cM=S^{d}$ in $\bbR^{d+1}$ for $d=5, 10, 15, 20$. \autoref{tab:simTime} shows that the RFPCA scales well for larger dimensions in terms of running time, and its relative disadvantage in speed as compared to the $L^2$~FPCA becomes smaller as $d$ and $n$ get larger.

The samples were generated in the same fashion as in the main text, except for the mean function $\muM(t) = \exp_{p_0}(2(d-1)^{-1/2}t, \dots, 2(d-1)^{-1/2}t, 0.3\pi\sin(\pi t), 0)$, and eigenfunctions $\phi_k(t) = d^{-1/2} R_t [\zeta_k(t / d), \dots, \zeta_k(t/d + (d-1)/d), 0]^T$, where $p_0 = [0, \dots, 0, 1]$ and $R_t$ is the rotation matrix from $p_0$ to $\muMt$. 

\single 
\begin{table}[ht]
\centering
\caption{A comparison of mean running time for $S^d$. The standard errors are below 2\% of the means.} 
\label{tab:simTime}
\begin{tabular}{c|cccc|cccc|cccc|cccc}
	           & \multicolumn{4}{|c|}{$n=50$} & \multicolumn{4}{|c}{$n=100$} & \multicolumn{4}{|c}{$n=200$} & \multicolumn{4}{|c}{$n=400$} \\
	   $d$    &  5   &  10  &  15  &   20    &  5   &  10  &  15  &   20    &  5   &  10  &  15  &   20    &  5   &  10  &  15   &   20   \\ \hline
	  RFPCA    & 1.3 & 1.7 & 2.1 &  2.7   & 1.9 & 2.6 & 3.3 &  4.3   & 3.1 & 4.6 & 5.7 &  7.4   & 5.9 & 7.8 & 10.5 & 13.2  \\
	$L^2$~FPCA & 0.4 & 0.7 & 1.0 &  1.5   & 0.8 & 1.2 & 1.8 &  2.4   & 1.4 & 2.5 & 3.5 &  4.4   & 3.0 & 4.4 & 6.3  &  8.2
\end{tabular}
\end{table}
\double 

\end{document}